\documentclass[a4paper,10pt]{article}

\usepackage{mathrsfs} 

\usepackage{amsthm}
\usepackage{amssymb}
\usepackage{latexsym}
\usepackage{amsmath}
\usepackage{amsfonts}
\usepackage{mathabx}
\usepackage[dvips]{graphicx}
\usepackage[utf8]{inputenc}

\usepackage[LGR,T1]{fontenc}%
\usepackage[francais,english]{babel}
\usepackage{url}
\usepackage{enumerate}
\usepackage{color}
\usepackage[urlcolor=cyan, colorlinks=true,  linkcolor=blue, citecolor=green]{hyperref}	
\usepackage{cleveref}
\usepackage{float}

\usepackage{fancyhdr}
\usepackage{tikz}
\usetikzlibrary{arrows.meta, calc,  fit,  positioning,  quotes}

\makeatletter
\tikzset{
  meta box/.style={
    draw,
    black,
    very thick,
    text centered
  },
  punkt/.style={
    meta box,
    rectangle,
    rounded corners,
    inner sep=3pt,
    minimum height=1em,
    minimum width=1em,
    align=center,
  },
  round box/.style={
    meta box,
    circle
  },
  every fit/.style={
    draw,
    thick,
    dashed,
    gray,
    inner sep=9pt
  }
}

\newcommand\tikz@expand@dimen[2]{\tikzset{minimum #2=#1}}
\tikzset{
  add dimen/.code 2 args={%
    \pgfkeysgetvalue{/pgf/minimum #1}\tikz@dimen@min
    \expandafter\tikz@expand@dimen\expandafter{\tikz@dimen@min + #2 * 2em}{#1}%
  },
  wider/.style={add dimen={width}{#1}},
  higher/.style={add dimen={height}{#1}},
}
\makeatother

\newtheorem{Theorem}{Theorem}[part]
\newtheorem{Definition}{Definition}[part]
\newtheorem{Proposition}{Proposition}[part]

\newtheorem{Assumption}{Assumption}[part]
\newtheorem{Lemma}{Lemma}[part]
\newtheorem{Corollary}{Corollary}[part]
\newtheorem{Remark}{Remark}[part]

\makeatletter \@addtoreset{equation}{section}

\@addtoreset{Definition}{section}

\@addtoreset{Theorem}{section}

\@addtoreset{Proposition}{section}

\@addtoreset{Property}{section}

\@addtoreset{Assumption}{section}

\@addtoreset{Corollary}{section}

\@addtoreset{Lemma}{section}

\@addtoreset{Remark}{section}

\@addtoreset{Example}{section}

\usepackage{bbm}
\usepackage{diagbox}
\usepackage{multicol}
\usepackage{multirow}
\usepackage{relsize}
\usepackage{tabularx}

\newcommand{\cA}{\mathcal{A}}

\newcommand{\cC}{\mathcal{C}}

\newcommand{\cE}{\mathcal{E}}
\newcommand{\cF}{\mathcal{F}}
\newcommand{\cG}{\mathcal{G}}
\newcommand{\cH}{\mathcal{H}}

\newcommand{\cL}{\mathcal{L}}
\newcommand{\cM}{\mathcal{M}}
\newcommand{\cN}{\mathcal{N}}

\newcommand{\cS}{\mathcal{S}}
\newcommand{\cT}{\mathcal{T}}
\newcommand{\cU}{\mathcal{U}}
\newcommand{\cV}{\mathcal{V}}
\newcommand{\cW}{\mathcal{W}}
\newcommand{\cX}{\mathcal{X}}
\newcommand{\cY}{\mathcal{Y}}
\newcommand{\cZ}{\mathcal{Z}}

\renewcommand{\P}{\mathbb{P}}

\newcommand{\R}{\mathbb{R}}

\def \proof{{\noindent \bf Proof. }}
\def \eproof{\hbox{ }\hfill$\Box$}

\newcommand{\ud}{\mathrm{d}}

\newcommand{\HYP}[1]
    {\ensuremath{({H#1} ) }}

\newcommand{\set}[1]
    {\ensuremath{\{ #1 \}}}
    
\newcommand{\HP}[1] 
    {\ensuremath{\mathscr{H}^{#1}}}

\newcommand{\esp}[1]{\ensuremath{\mathbb{E} \!\! \left[#1\right] }}

\newcommand{\EFp}[2]
    {\ensuremath{
     \mathbb{E}_{#1}\!\!\left[#2\right] }}
     
\renewcommand{\Xi}[1]{X_{i #1}}

\definecolor{vert}{HTML}{000000}

\newcommand{\mfy}{\mathfrak{y}}
\newcommand{\mfz}{\mathfrak{z}}
\newcommand{\mfa}{\mathfrak{a}}
\newcommand{\mfc}{\mathfrak{c}}

\title{
Deep Runge-Kutta schemes for BSDEs
}

\author{Jean-François Chassagneux \thanks{Universit\'e Paris Cit\'e, Laboratoire de Probabilit\'es, Statistiques et Mod\'elisation, F-75013 Paris, France. Email: chassagneux@lpsm.paris.} \, ,  Junchao Chen  \thanks{Universit\'e Paris Cit\'e, Laboratoire de Probabilit\'es, Statistiques et Mod\'elisation, F-75013 Paris, France.  Email: juchen@lpsm.paris} \,, Noufel Frikha \thanks{Universit\'e Paris 1 Panth\'eon-Sorbonne, Centre d'Economie de la Sorbonne, 106 Boulevard de l'H\^opital, 75642 Paris Cedex 13, France, Email: noufel.frikha@univ-paris1.fr}.}

\date{}

\begin{document}
\maketitle

\begin{abstract}
{\color{vert}
\textcolor{black}{We propose a new probabilistic scheme which combines deep learning techniques with high order schemes for backward stochastic differential equations belonging to the class of Runge-Kutta methods to solve high-dimensional semi-linear parabolic partial differential equations.} Our approach notably extends the one introduced in \cite{hure2020deep} for the implicit Euler scheme to schemes which are more efficient in terms of discrete-time error. We establish some convergence results for our implemented schemes under classical regularity assumptions. We also illustrate the efficiency of our method for different schemes of order one, two and three. Our numerical results indicate that the Crank-Nicolson schemes is a good compromise in terms of precision, computational cost and numerical implementation.
%
}
\end{abstract}


\section{Introduction}
In this paper, we study the numerical approximation of Backward Stochastic Differential Equations (BSDEs for short) or equivalently semi-linear parabolic Partial Differential Equations (PDEs) by combining high order discretization schemes with deep learning techniques. We work in a Markovian setting introducing  a forward diffusion process $(\cX)$ with dynamics 
\begin{equation}
 \cX_t  = \cX_0 + \int_0^t \mu(\cX_s) \ud s +   \int_0^t  \sigma(\cX_s) \ud W_s, \,\quad 0 \le t  \le T, \label{eq forward sde chapter 3}
\end{equation}

\noindent and we aim to approximate the solution $(\cY,\cZ)$ of the BSDE
\begin{align}
\cY_t & = g(\cX_T) + \int_t^Tf(\cX_s,\cY_s,\cZ_s) \ud s - \int_t^T \cZ_s \cdot \ud W_s\,,\;0\le t \le T\,,
\label{eq backward sde chapter 3}
\end{align} 
where $ W$ is a $d$-dimensional Brownian motion defined on a complete probability space $(\Omega , \cA,\P)$, $\mu : \R^d \rightarrow \R^d$ and $\sigma: \R^d \rightarrow \mathbb{M}_d$ (the set of $d\times d$ matrices) are measurable functions and the initial condition $\cX_0$ is an $\mathbb{R}^d $-valued  random variable \textcolor{black}{independent of $W$}. We denote  the  filtration generated by $W$ and $\cX_0$ as $(\mathcal{F}_{t})_{0\leq t \leq T}$, augmented with $\P$ null sets. 

Our approach relies on the classical connection between (BSDEs) and semi-linear parabolic PDEs as initiated in \cite{Pardoux:Peng:90}. Namely, under some standard regularity assumptions on $\mu$, $\sigma$, $f$ and $g$, it holds
\begin{align}
	Y_t & = u(t, \cX_t),  \quad Z_t = \sigma^\top(\cX_t) \nabla_x u(t, \cX_t) , \quad 0 \leq t \leq T,
\end{align}
where $u: [0, T] \times \mathbb{R}^d \longmapsto \mathbb{R} $ is the solution to the following parabolic semi-linear PDE:
\begin{equation}
\label{eq pde chapter 3}
\left\{
\begin{array}{ll}
\partial_t u(t, x) + \mathcal{L} u(t, x) + f(u(t,x), \sigma^{\top}(x)\nabla_x u(t, x))  =0, &   (t, x) \in [0,T) \times \mathbb{R}^d, \\
 u(T, x) = g(x), \quad x \in \mathbb{R}^d  & 
\end{array}
\right.
\end{equation}
and $\mathcal{L}$ is the infinitesimal generator defined by
\begin{equation}
\label{infinitesimal:operator:diffusion}
\mathcal{L} u(t, x) :=  \mu(x)\cdot \nabla_x u(t, x) +\frac12 \mathrm{Tr}[(\sigma\sigma^{\top})(x)\nabla^2_x u(t, x)].
\end{equation}

Since BSDEs have been introduced by Pardoux and Peng \cite{Pardoux:Peng:90, Pardoux:Peng:92} in 1990s, their numerical approximation has attracted considerable attention. However, solving high-dimensional BSDEs is still a challenging task. A first stream of methods relies on the backward programming algorithm (the method we study in this work belongs to this class). The main difficulty for these methods is  to compute the conditional expectations. Many solutions have been proposed in the last two decades, such as the cubature methods \cite{chassagneux2020cubature, crisan2012solving, crisan2014second}, optimal quantization methods \cite{bally2003, BALLY20031, PAGES2018847, nmeir2021quantization}, Malliavin calculus based methods \cite{CRISAN, BOUCHARD2004175, hu2011} and some linear regression methods \cite{gobet2005regression, gobet2016stratified, gobet2016approximation}. These methods are efficient in a low dimensional setting $d \sim 5$. However they all face the ``curse of dimensionality''. Other type of methods have also been recently introduced, to quote a few: the branching process method (see e.g. \cite{henry2014numerical}), the forward method using Wiener chaos expansion, that allows to handle possibly  non Markovian setting (see e.g. \cite{briand2014simulation}) and the Multi-Level Picard (MLP) method (see e.g. \cite{hutzenthaler2019multilevel}) which does not suffer from the curse of dimensionality.
%

In the past five years, numerical procedures using  learning methods (linear and non-linear regression) have been proposed with excellent empirical results even in high dimensional context. This includes  forward schemes \cite{weinan2017deep, beck2017machine, han2020convergence, chassagneux2021learning} or  backward schemes \cite{hure2020deep, germain2020deep, pham2021neural,germain2021neural}, see also \cite{chan2019machine, jiang2021convergence, kapllani2020deep, takahashi2021new, chassagneux2021numerical}. However, these algorithms rely on the standard Euler-Maruyama time discretization scheme whose weak convergence rate is of order 1, so that the overall computational time cost for these algorithms can still be large for high-dimensional BSDEs as many time steps might be required to achieve good accuracy. 

In this paper, we combine some high-order time discretization numerical schemes for BSDEs with non-linear regression based on (deep) neural networks to solve high-dimensional \textcolor{black}{semi-linear parabolic PDEs}. Many high-order discrete-time approximation schemes have been introduced for BSDEs see among others \cite{zhao2006new, zhao2010stable, chassagneux2014linear,chassagneux2014runge}.
We here chose to work with the family of Runge-Kutta methods \cite{chassagneux2014runge} which are one step schemes. It includes the implicit Euler scheme, explicit Euler scheme and the Crank-Nicolson scheme, but also schemes with more than one stage of computation per step, see Definition \ref{de RK scheme} below. 
These high-order schemes are based on a backward algorithm and, as usual, they require a good approximation procedure of conditional expectations in practice. We chose here to compute them using deep learning techniques as previous papers demonstrated their efficiency to tame the curse of dimensionality.
To the best of our knowledge, high order schemes for BSDEs have not been tested with (non-linear) regression techniques.
Our main contribution is thus to present an approach to implement the Runge-Kutta schemes using non-linear regression techniques in an efficient manner. Our method is inspired by the one introduced in \cite{hure2020deep} for the implicit Euler scheme. We also establish an upper-bound for the global error of our scheme which writes as the sum of the discrete time error and the approximation error induced by the use of neural networks. This is obtained by a careful analysis of the scheme, in particular its stability properties. We thus identify order one, two and three schemes whose efficiency is compared numerically. Based on our numerical experiments, we would recommend the use of the Crank-Nicolson scheme for its tractability, good precision and reasonable computational cost.


The rest of the paper is organized as follows. We first recall the definition of Runge-Kutta schemes for BSDEs in \Cref{se Runge-Kutta schemes}, then we study their stability in two different ways. Theorem \ref{th disc error control} provides the discrete time error for the main methods that will be studied in this paper. In \Cref{se a learning method}, we present an implementation of the Runge-Kutta schemes to solve BSDEs by neural networks, including the special case of the implicit Euler schemes \cite{hure2020deep}, explicit Euler scheme, Crank-Nicolson scheme, two stage and three stage explicit Runge-Kutta scheme. We also establish the global error bound for a general learning method applied to Runge-Kutta schemes, see Theorem \ref{th control global error}. In \Cref{se numerics}, we numerically illustrate the convergence order of the discrete time error of the methods presented in Section \ref{subse implemented schemes}. We also compare the computational time cost of these methods.
\vspace{1\baselineskip}

We conclude this introduction by presenting some standing assumptions and notations that will be used in the sequel.
From now on, we assume that the driver  $f$, the terminal function $g$ and the coefficients $\mu$ and $\sigma$ of the SDE \eqref{eq forward sde chapter 3} satisfy the following regularity assumption:
\begin{Assumption}
	There exists $0 \le K< \infty$ such that for all $x_1, x_2 , z_1, z_2 \in \mathbb{R}^d$, $y_1, y_2 \in \mathbb{R}$
	\begin{align}
		|f(x_2, y_2, z_2) - f(x_1, y_1, z_1)| & \le K \left( |x_1-x_2| + |y_1-y_2| + |z_1 - z_2| \right), 
		\label{eq lip f}
		\\
		|g(x_2) - g(x_1) | & \le K|x_1 -x_2 |,
	\end{align}
	and
	\begin{align}
	&	| \mu( x_1)  -  \mu( x_2) | + | \sigma( x_1)  -  \sigma(x_2) | \le K |x_1-x_2|. 
\end{align}
\end{Assumption}

%
%



\subsection{Notations related to neural networks}
\label{subse NN notations}
We now introduce some notations and basic definitions concerning the class of neural networks that will be used in this paper:
\begin{align*}
	 \left\{
	 		\begin{array}{rl}
	 			\mathbb{M}_{p\times q}(\mathbb{R}): &\text{the matrix space for all } p \times q \text{ matrices with elements in } \mathbb{R} ,\\
	 			d_0 = d : &\text{input dimension,} \\
	 			d_1: &\text{output dimension,}\\
	 			L + 1 \in \mathbb{N} \setminus \{0,1,2\}: & \text{number of layers  of the network,} \\
	 			m_{\ell},  \ell = 0, 1,  \cdots, L: & \text{number of neurons on each layer, note } m_0 = d_0, m_L = d_1,\\
	 			\cH_n^\ell, \ell = 1,  \cdots, L-1: & \text{output of the hidden layers at time $t_n, 0\le n \le N-1$}.
	 		\end{array}
	 \right.
\end{align*}

For the $L-1$ hidden layers in this neural network, we choose for simplicity the same number of neurons $m_\ell = m, \ell = 1, \cdots, L-1$. \\

For $ \ell = 1, \cdots, L $, we define the maps $\cM_\ell : \mathbb{R}^{m_{\ell - 1}} \longmapsto \mathbb{R}^{m_{\ell}}  $ as:
\begin{align}
	\cM_\ell (x) = \cW_\ell x + \beta_\ell,
\end{align}
where $ \cW_\ell \in \mathbb{M}_{m_\ell \times m_{\ell-1}}(\mathbb{R} ) $ is a matrix called weight, and $\beta_\ell \in \mathbb{R}^{m_\ell} $ is a vector called bias. Then $ \cM_\ell $ is an affine transformation that can map the features of the $(\ell - 1)\verb|-|th$ layer to the $\ell \verb|-|th$ layer.
A feedforward neural network is a function from $\mathbb{R}^{d_0}$ to $\mathbb{R}^{d_1}$ defined as the composition
\begin{align}\label{eq de NN function}
	\mathbb{R}^{d_0} \ni x \longmapsto	\cN_m (x; \theta) = \cM_L \circ \rho_{L-1} \circ \cM_{L-1} \circ \cdots \circ \rho_1 \circ \cM_{1}(x) \in \mathbb{R}^{d_1},
\end{align}

\noindent where $\rho_\ell (x) = ( \rho(x_1), \cdots, \rho(x_{m_\ell}) ), x\in \mathbb{R}^{m_\ell}, \ell = 1, \cdots, L-1 $, here $ \rho: \mathbb{R} \longmapsto \mathbb{R} $ is an activation function which is also a nonlinear function, such as ReLu, Elu, tanh, sigmoid. Then the parameters of the neural network consist of the weight matrices $ (\cW_\ell)_{1\le \ell \le L} $, the bias vector $ (\beta_\ell)_{1\le \ell \le L} $. For fixed $d_0, d_1$ and $L$, the total number of parameters is 
\[ N_m := \sum\limits_{\ell = 1}^{L}  m_\ell (m_{\ell-1} + 1 ) = d_0(1+m) + m(m+1)(L-2) + m(1+d_1) ,\]
 so that the parameters can be identified with an element $ \theta \in \mathbb{R}^{N_m} $. We will sometimes insist on the output dimension of the neural network by using $N_m^{d_1=q}$, for a $q$ dimensional output. In particular, let us mention that for Euler schemes, the output dimension is $d_1 := 1+d$ which includes one component for the $Y-$part and $d$ components for the $Z-$part. But for general Runge-Kutta schemes, we will set $d_1 := 1+2d$ for the networks: The components consist  $1, d, d$ dimensions for $Y, Z, A$, respectively, where $A$ is an extra $d-$dimensional variable that will be required to define the scheme, see Section \ref{subse implementation RK schemes}. 

We introduce 
\begin{align}\label{eq set of networks}
	 \cS_{d_0, d_1, L, m}^{\rho}(\mathbb{R}^{N_m})  := \big\{ \cN_m(\cdot; \theta):\R^d \rightarrow  \mathbb{R}^{d_1} \quad | \quad  \theta \in \mathbb{R}^{N_m} \big\} 
\end{align}
and 
\begin{align}
	\cS_{d_0, d_1, L}^{\rho} := \bigcup_{m \in \mathbb{N}^+} \cS_{d_0, d_1, L, m}^{\rho}(\mathbb{R}^{N_m}).
\end{align}

\noindent
We conclude this introductory section by recalling the fundamental result of Hornik et al. \cite{hornik1989multilayer} that states the following universal approximation theorem and  justifies that the neural networks can be applied as function approximators:

\begin{Theorem}[Universal approximation theorem] \label{universal approximation theorem}
	$\cS_{d_0, d_1, L}^{\rho}$ is dense in $L^2(\upsilon )$ for any finite measure $\upsilon$ on $ \mathbb{R}^d $, whenever $\rho$ is continuous and non-constant.
\end{Theorem}


\section{Runge-Kutta schemes for BSDEs} \label{se Runge-Kutta schemes}


{\color{vert}
We work with the class of Runge-Kutta schemes that have been introduced in \cite{chassagneux2014runge} in the BSDE setting. The main difference with our approach is that we also consider an approximation of the forward process. In this section, we first recall the definition of the Runge-Kutta schemes for BSDEs. We then study the stability properties of this class of scheme. We conclude the section by analyzing the discrete time error of the main schemes of interest.
}

\subsection{Definitions}
For ease of presentation, we consider an equidistant time grid 
$$ \pi := \{ t_0 = 0 < \cdots < t_n < \cdots < t_N = T \} $$ 

\noindent of the interval $[0,T]$ with time step $h: = \frac{T}{N}, t_n = nh, n= 0, \cdots, N$ and let $ \Delta W_n = W_{t_{n+1}} - W_{t_n} , 0 \le n \le N-1$.

\vspace{2mm}
The Runge-Kutta schemes involve in full generality intermediate steps of computation between two dates of the main grid $\pi$.
Thus, for a positive integer $Q$, let $ c = (c_1, \ldots, c_{Q+1}) \in [0, 1]^{Q+1} $ satisfying $ 0 =: c_1 < c_2 \le \ldots \le c_q \le \cdots \le c_Q \le c_{Q+1} := 1 $. We introduce the intermediate ``instances'' $ t_{n, q} : = t_{n+1} - c_q h $, for $0\leq n\leq N-1$. With these notations, we observe that $ t_n = t_{n, Q+1} \le \ldots \le  t_{n,q} \le \ldots \le t_{n, 1} = t_{n+1} $. We denote the ``full grid'' 
$$\Pi := \set{ t_{n,q} \in [0,T]\,|\, 0 \le n \le N-1, 1 \le q \le Q+1}.$$

First, we are given an approximation of the forward component \eqref{eq forward sde chapter 3} on the grid $\Pi$. Namely, for $t_{n,q} \in \Pi$, $\cX_{t_{n,q}}$ is approximated by $X_{n,q} \in \cL^2(\cF_{t_{n,q}})$, $0 \le n \le N-1$ and $ 1 \le q \le Q+1 $. For ease of notation, we will simply denote by $(X_n)_{0 \le n \le N}$ the approximation of $\cX$ on the grid $\pi$. Observe that $X_{n,Q+1}=X_n$ and $X_{n,1}=X_{n+1}$, for $0 \le n \le N-1$. In the following, we assume that $X$ is a Markov process on $\Pi$. 
We denote $ \Delta W_{n, q} = W_{t_{n+1}} - \textcolor{black}{ W_{t_{n, q}}} , 0 \le n \le N-1$, $1 \le q \le Q+1$.

\vspace{2mm}
\noindent We now define $(Y,Z)$ the approximation of $(\cY,\cZ)$ given by \eqref{eq backward sde chapter 3}.
\begin{Definition}\label{de RK scheme}
\begin{enumerate}[i)]
	\item \textcolor{black}{Assuming that $g$ is differentiable,} we set the terminal condition to  $$ (Y_N, Z_N) = (g(X_N), \sigma(X_N)^\top \nabla_x g(X_N) ) .$$
	\item For $0 \le n \le N-1$ and $Q\ge 1$, the transition from $(Y_{n+1}, Z_{n+1})$ to $(Y_n, Z_n)$ involves $Q$ stages.
	At the intermediate instances, for $1 < q \le Q+1$, let
	 \begin{align}
	 	Y_{n,q} &= \EFp{t_{n,q}}{ Y_{n+1} + h \sum_{k=1}^{q} a_{qk} f(X_{n,k},Y_{n,k}, Z_{n,k}) },
		\label{eq de interm stage Y}
		\\
	 	Z_{n,q} &= \EFp{t_{n,q}}{H_q^n Y_{n+1} + h \sum_{k=1}^{q-1} \alpha_{qk} H_{q,k}^n f(X_{n,k},Y_{n,k}, Z_{n,k}) },
		\label{eq de interm stage Z}
	 \end{align}
	 where $ (a_{qk})_{1\le q, k\le Q+1}, (\alpha_{qk})_{1\le q, k\le Q+1}$
	  take their values in $\R$ and
	  with $ a_{1k} = \alpha_{1k} = 0, 1\le k \le Q  $, $ a_{qk} = \alpha_{qk} = 0 ,  1\le q < k \le Q+1$ and
	 \begin{align}
	 	\sum_{k=1}^q a_{qk} = \sum_{k=1}^{q-1} \alpha_{qk} \mathbbm{1}_{\{c_k < c_q\}} = c_q, \qquad q \le Q+1. 
	 \end{align}
	We then set $(Y_n,Z_n) = (Y_{n,Q+1},Y_{n,Q+1})$, \textcolor{black}{for $0\leq n \leq N-1$}.\\
	 For any $1\le k < q\le Q+1$, $0\leq n \le N-1$, the random variables $H_q^n, H_{q,k}^n $ are $ \cF_{t_{n+1}} -$measurable
	 and satisfy:
	 \begin{align}
	 	&\EFp{t_{n,q}}{H_q^n} = \EFp{t_{n,k}}{H_{q,k}^n} = 0 \text{ and } \upsilon^n_q := \EFp{t_{n,q}}{|H_q^n|^2}, 
		 \upsilon^n_{q,k} := \EFp{t_{n,k}}{|H_{q,k}^n|^2},   
		 \label{eq conditions on H 1}
		 \\
		 & \EFp{t_{n,q}}{H_q^n \Delta W_{n, q}} = 1,
		  \label{eq conditions on H 3}
		   \\
		 & \frac{\lambda}h \le \min(\upsilon^n_q, \upsilon^n_{q,k} ) \text{ and } \max(\upsilon^n_q, \upsilon^n_{q,k} ) \le \frac{\Lambda}{h},
		  \label{eq conditions on H 2}
	 \end{align}
	 where $\lambda, \Lambda $ are positive constants which do not depend on $h$. \\
\end{enumerate}
\end{Definition}

%

{\color{vert}
\noindent We observe that \eqref{eq conditions on H 3} is mainly a renormalisation property that will ease the presentation of our method. We note also that \eqref{eq de interm stage Y} may define $Y_{n,q}$ implicitly. For the rest of the paper, we assume
\begin{align}\label{eq cond well posedness of the scheme}
\max\limits_{1<q<Q+1}|a_{qq}|hK<1,
\end{align}
which guarantees that the scheme is well-defined, where $K$ is defined by \eqref{eq lip f}. In particular, by a direct induction argument, it holds
\begin{align}
\max_{n,q}\esp{|Y_{n,q}|^2 + |Z_{n,q}|^2} < +\infty \,.
\end{align}
}


%

\subsection{Stability of Runge-Kutta scheme}
A key property to obtain the convergence results stated in Theorem \ref{th disc error control} is -- classically -- the $L^2$-stability of the schemes of Definition \ref{de RK scheme}. This has already been observed in \cite{chassagneux2014runge}. We shall review here this property as it will be useful in the sequel.

\vspace{2mm} The first observation is the fact that the schemes given in Definition \ref{de RK scheme} can be written in the following implicit form, for $0\leq n \leq N-1$:
\begin{align}
Y_n &= \EFp{t_n}{Y_{n+1} + \Phi^Y_n(Y_{n+1},Z_{n+1},h)},
\\
Z_n &= \EFp{t_n}{H_{Q+1}^nY_{n+1} + \Phi^Z_n(Y_{n+1},Z_{n+1},h)},
\end{align}

\noindent where $(\Phi^Y_n,\Phi^Z_n): \Omega \times \cL^2(\cF_{t_{n+1}})\times \cL^2(\cF_{t_{n+1}}) \times \R_{+} \rightarrow \cL^2(\cF_{t_{n+1}})$.
This writing really stresses the fact that the schemes are one-step scheme. 
One introduces a perturbed version of the scheme, namely,
\begin{align}
\check{Y}_n &= \EFp{t_n}{\check{Y}_{n+1} + \Phi^Y_n(\check{Y}_{n+1},\check{Z}_{n+1},h)} + \check{\zeta}^Y_n
\\
\check{Z}_n &= \EFp{t_n}{H_{Q+1}^n\check{Z}_{n+1} + \Phi^Z_n(\check{Y}_{n+1},\check{Z}_{n+1},h)} + \check{\zeta}^Z_n
\end{align}

\noindent for $(\check{\zeta}^Y_n,\check{\zeta}^Z_n) \in \cL^2(\cF_{t_n})\times\cL^2(\cF_{t_n})$, 
and obtains, using Theorem 1.2(i) in \cite{chassagneux2014runge}, the following stability result, setting $\delta Y_n :=\check{Y}_n - Y_n$, $\delta Z_n :=\check{Z}_n - Z_n$,
\begin{align}\label{eq stab for disc time error}
\max_{n<N}{\esp{|\delta Y_n|^2}} + \sum_{n=0}^{N-1}h\esp{|\delta Z_n|^2}
\le
C\esp{|\delta Y_N|^2 + h{|\delta Z_N|^2} + \sum_{n=0}^{N-1}\frac{|\check{\zeta}^Y_n|^2}{h} + h|\check{\zeta}^Z_n|^2}\,.
\end{align}
This approach is particularly well-suited for the study of the discrete time error, see the proof of Theorem \ref{th disc error control} in Section \ref{subse discrete time error} below.
However, we need also a stability result to control the error linked to the estimation of the conditional expectations at each stage of the schemes. To this end, we now introduce another perturbed scheme, for $0\leq n \leq N-1$, at the intermediate instances, for $1 < q \le Q+1$, let
	 \begin{align}
	 	\tilde{Y}_{n,q} &= \EFp{t_{n,q}}{ \tilde{Y}_{n+1} + h \sum_{k=1}^{q} a_{qk} f(X_{n,k},\tilde{Y}_{n,k}, \tilde{Z}_{n,k}) } + \zeta^y_{n,q},
		\label{eq de interm stage Y perturbed}
		\\
	 	\tilde{Z}_{n,q} &= \EFp{t_{n,q}}{H_q^n \tilde{Y}_{n+1} + h \sum_{k=1}^{q-1} \alpha_{qk} H_{q,k}^n f(X_{n,k},\tilde{Y}_{n,k}, \tilde{Z}_{n,k}) } + \zeta^z_{n,q},
		\label{eq de interm stage Z perturbed}
	 \end{align}
	 with $(\zeta^y_{n,q},\zeta^z_{n,q}) \in \cL^2(\cF_{t_{n,q}})$.

\vspace{2mm}
\noindent Associated to the above perturbed version, we can state the following stability result, 
whose proof is postponed to the Appendix, see Section \ref{subse proof of prop}.
\begin{Proposition}\label{pr generic stab for cond exp estim} Let $\delta Y_n :=\tilde{Y}_n - Y_n$, $\delta Z_n :=\tilde{Z}_n - Z_n$. There exists $0 \le C< \infty$ such that for $h$ small enough, it holds
\begin{align}\label{eq stab for cond expec estim}
\max_{0\leq n \leq N-1}{\esp{|\delta Y_n|^2}} + \sum_{n=0}^{N-1}h\esp{|\delta Z_n|^2}
\le
C\esp{|\delta Y_N|^2 + h{|\delta Z_N|^2} + 
\sum_{n=0}^{N-1}\sum_{q=2}^{Q+1} \left(\frac{|\zeta^y_{n,q}|^2}h + h|\zeta^z_{n,q}|^2\right) }\,.
\end{align}
\end{Proposition}
%

\subsection{Discrete time error} \label{subse discrete time error}

In \cite{chassagneux2014runge}, the discrete-time error has been studied when $X=\cX$, namely there is no error in the approximation of the underlying process. Building on the results in \cite{chassagneux2014runge}, we will give an upper bound of the discrete-time error when the forward process is itself approximated. We will focus here on one stage schemes both implicit and explicit and two and three stage explicit schemes. We refer to Remark \ref{re restriction} for more details on this limitation. 
The control of the discrete-time error is based on smoothness assumptions satisfied by the value function $u$ solution to \eqref{eq pde chapter 3}. We now introduce the necessary notations to formalise our statement.

\noindent Let 
$$\cM := \set{\emptyset} \cup \bigcup_{m=1}^\infty \set{0,\dots,d}^m, $$
the set of multi-indices with entry $0, \dots, d$.
We define the  differential operators as 
\begin{align}
	\cL^{(0)} &:= \partial_t + \sum_{i=1}^d \mu_i \partial_{x_i} + \frac12 \sum_{i=1}^d \sum_{j=1}^d (\sigma \sigma^\top)_{ij} \partial_{x_i, x_j}^2 = \partial_t + \cL , \\
	\cL^{(\ell)} &:= \sum_{k=1}^d \sigma_{k \ell } \partial_{x_k}, \qquad \ell \in \{1, \ldots, d\},
\end{align}
and their iteration, namely, for $\alpha \in \cM$,
$$
\cL^\alpha := \cL^{(\alpha_1)}\circ\dots\circ \cL^{(\alpha_p)},
$$
for a multi-index $\alpha$ with positive length $\mathfrak{l}(\alpha):=p$. By convention, $L^{\emptyset}$ is the identity operator ($\mathfrak{l}(\emptyset)=0$), and denote $ \alpha = (0, \cdots, 0) =: (0)_p $ with positive length $\mathfrak{l}(\alpha) =p$. By convention, we set $(0)_0:=\emptyset$. We denote by $*$ the concatenation of two multi-indices namely $\alpha*\beta = (\alpha_1,\dots,\alpha_p,\beta_1,\dots,\beta_q)$ with positive length $p = \mathfrak{l}(\alpha)$ and positive length $q= \mathfrak{l}(\beta)$. We set $\emptyset * \alpha = \alpha * \emptyset = \alpha$, for all $\alpha \in \cM$.
We denote by $ \cG_b^{l}$ the set of all functions $ v:[0, T] \times \mathbb{R}^d \longmapsto \mathbb{R} $ for which $\cL^\alpha v$ is well defined, continuous and bounded for all muti-index $\alpha \in  \{ (\alpha_1, \cdots, \alpha_p) | 1\le p \le l\} $. We also denote classicaly $\cC^p_b$ the set of $p$ times continuously differentiable functions with all their derivatives bounded.

In particular, we shall use the following assumption, for $p = 1,2,3$:\\
\\ $(\mathrm{H}r)_p:$ \textcolor{black}{The solution $u$ to the PDE \eqref{eq pde chapter 3} belongs to $ \cG_b^{p+1}$} and $f \in \cC^p_b$. 

\vspace{2mm} The key to control the discrete-time error by a judicious choice of the scheme coefficients is to be able to expand the value function $u$ along the approximation scheme $X$. To this end, we introduce the following  assumption for a positive integer $M$:\\
\\ $(\mathrm{H}X)_M:$ the process $X$ satisfies, for all $v \in \cG^{M+1}_b$, $0 \le n \le N-1$, $1 < q \le Q+1$, $k \le q$, $1\le \ell \le d$, denote { $ v^\alpha = \cL^\alpha v $},
\begin{align}
&\EFp{t_{n,q}}{v(t_{n,k},X_{n,k})} = \sum_{m=0}^M v^{(0)_m}(t_{n,q}, X_{n,q})\frac{(\set{c_q-c_k}h)^m}{m!} + O_{{t_{n,q}}}(h^{M+1}),
\label{eq approx exp y}
\\
&\EFp{t_{n,q}}{(H_{q}^n)^\ell v(t_{n+1},X_{n+1})} = \sum_{m=0}^{M-1} v^{(\ell)*(0)_m}(t_{n,q}, X_{n,q})\frac{(c_qh)^m}{m!} + O_{{t_{n,q}}}(h^{M}),
\label{eq approx exp z 1}
\\
&\EFp{t_{n,q}}{(H^n_{q,k})^\ell v(t_{n+1},X_{n+1})} = \sum_{m=0}^{M-1} v^{(\ell)*(0)_m}(t_{n, q}, X_{n,q})\frac{(\set{c_q-c_k}h)^m}{m!}  + O_{{t_{n,q}}}(h^{M}),
\label{eq approx exp z 2}
\end{align}
 where  the notation $R=O_{t}(r)$, $r > 0$, means that the random variable $R$ is such that $ |R| \le \lambda_t^r r $ with  $ \lambda_t^r  $ is a positive random variable satisfying for all $p, r>0$
\begin{align*}
	\esp{|\lambda_t^r |^p} \le C_p, 
\end{align*}
\noindent \textcolor{black}{for some constant $C_p <\infty$ independent of $r$.}
We note that \eqref{eq approx exp y} indicates that $X$ is a weak approximation scheme of order $M$ approximation \cite{kloeden1992stochastic, chassagneux2014runge}. Conditions \eqref{eq approx exp z 1}-\eqref{eq approx exp z 2} are required to handle the error coming from the approximation of the $Z$-component.


\noindent Regarding the discrete-time error, our main result reads as follows.
Define 
\begin{align}\label{eq de support bsde} 
(\bar{Y}_n,\bar{Z}_n) = (u(t_n,X_n), \sigma^\top \nabla_x u(t_n,X_n)) ,
\end{align}
and the global discrete time error as
\begin{align}\label{eq de truncation error}
\cT_N = \max_{0\leq n \leq N-1}\esp{|\bar{Y}_n - Y_n|^2} + h \sum_{n=0}^{N-1}\esp{|\bar{Z}_n - Z_n|^2} .
\end{align}
As usual, we say that the scheme is of order $\mathfrak{a} \in [0,\infty)$ if $\cT_N = O(h^{ 2 \mathfrak{a} })$.
\begin{Theorem}\label{th disc error control}
\begin{enumerate}[i)]
\item Theta-scheme (one stage scheme): Assume $(\mathrm{H}r)_1$ and $(\mathrm{H}X)_1$, the global \textcolor{black}{discrete time} error of Runge-Kutta scheme is at least of order 1 if 
	\begin{align*}
		a_{21} + a_{22} = 1.
	\end{align*} 
	This condition leads to the explicit Euler scheme when $ a_{21} = 1 $ and implicit Euler scheme when $a_{22} = 1$, respectively.\\
\item Crank-Nicolson scheme (one stage scheme): Assume $(\mathrm{H}r)_2$ and $(\mathrm{H}X)_{2}$, the global \textcolor{black}{discrete time} error of Runge-Kutta scheme is at least of order 2 if 
	\begin{align*}
		a_{21} = a_{22} = \frac{1}{2} \qquad and \qquad \alpha_{21} = 1.
	\end{align*} 
	
\item Two stage explicit scheme: Assume $(\mathrm{H}r)_2$ and $(\mathrm{H}X)_2$, the global \textcolor{black}{discrete time} error of Runge-Kutta scheme is at least of order 2 if $a_{22} = a_{33} = 0$ and
	\begin{align*}
		a_{21} = c_2, \quad a_{31} = 1 - \frac{1}{2c_2}, \quad a_{32} = \frac{1}{2c_2},  \quad and \quad \alpha_{31} + \alpha_{32}\mathbbm{1}_{\{c_2<1\}}=1.
	\end{align*}

\item Three stage explicit scheme: Assume $(\mathrm{H}r)_3$ and $(\mathrm{H}X)_3$, the global \textcolor{black}{discrete time} error of Runge-Kutta scheme is at least of order 3 if $0 < c_2 <c_3 \le 1 ( c_2 \neq \frac{2}{3}\; when \; c_3=1), a_{22} = a_{33} =a_{44}= 0$ and the following conditions holds true
	\begin{align*}
		& a_{41} + a_{42} + a_{43} = 1, \quad a_{42} c_2 + a_{43} c_3 = \frac{1}{2}, \\
		& a_{42}c_2^2 + a_{43}c_3^2 = \frac{1}{3}, \quad a_{43}a_{32}c_2 = a_{43} \alpha_{32}c_2 = \frac{1}{6}, \\
		& \alpha_{41} + \alpha_{42} + \alpha_{43} \mathbbm{1}_{\{ c_3 < 1\}} = 1, \quad \alpha_{42}  c_2 + \alpha_{43}c_3 \mathbbm{1}_{\{c_3 < 1\}} = \frac{1}{2}.
	\end{align*}

\end{enumerate}
\end{Theorem}
\proof
1. We give a short proof of the discrete time error upper bound for these most interesting schemes, as it is based on the one given in  \cite{chassagneux2014runge}.\\
Recall the definition of $(\bar{Y}_n,\bar{Z}_n)_{0 \le n \le N}$ in \eqref{eq de support bsde}. We observe that it can be written also as 
\begin{align*}
\bar{Y}_n &= \EFp{t_n}{\bar{Y}_{n+1} + \Phi^Y_n(\bar{Y}_{n+1},\bar{Z}_{n+1},h)} + \check{\zeta}^Y_n,
\\
\bar{Z}_n &= \EFp{t_n}{H_{Q+1}^n\bar{Y}_{n+1} + \Phi^Z_n(\bar{Y}_{n+1},\bar{Z}_{n+1},h)} + \check{\zeta}^Z_n.
\end{align*}
Now, thanks to assumption $(\mathrm{H}X)_M$ with $M=1, 2, 3$, one can follow the computations made in Theorem 1.3 for statement (i)-(ii), Theorem 1.5 for statement (iii) or Theorem 1.6 for statement (iv) in \cite{chassagneux2014runge} to obtain the corresponding upper bounds for the error linked to the perturbation, namely
\begin{align*}
\esp{\frac{|\check{\zeta}^Y_n|^2}{h} + h |\check{\zeta}^Z_n|^2} = O(h^{ 2 \mathfrak{a}+1}), 
\end{align*}
for $\mathfrak{a} = 1,2 \text{ or } 3$. The proof is then concluded using \eqref{eq stab for disc time error}. 
%
\eproof

\begin{Remark} \label{re restriction} 
	(i) As mentioned in \cite{chassagneux2014runge} and contrary to the ODE case, there exists an order barrier for implicit scheme to get an order $Q+1$ scheme with a $Q-$stage scheme when $Q>1$: This is indeed the case as soon as $\partial_z f \neq 0$, since the scheme given in Definition \ref{de RK scheme} are always explicit for the $Z$-component. Hence, we only consider the explicit scheme when $Q>1$ as the implicit scheme has no advantage compared to the explicit scheme for general drivers $f$.\\
	(ii) For the explicit Runge-Kutta scheme, we can choose the coefficients such that $ \alpha_{qk} = a_{qk}$, for all $ 1 \le k < q \le Q+1 $.\\
	(iii) We do not consider the case with $Q>3$ in this paper for two main reasons. First, there is an order barrier when $Q > 3$ (we need more than $Q$ stage to get an order $Q$ scheme), see  \cite{chassagneux2014runge} for a proof: This will make this scheme computationally very costly. Secondly, in our numerical experiments, when setting $Q=3$, the discrete time error is extremely small even for small number of time steps $N$ \textcolor{black}{so that the global error} is dominated by the variance. Hence, we cannot observe the convergence order of the scheme. We refer to Section \ref{se numerics} for more details on this specific point.
\end{Remark}

\section{A learning method for Runge-Kutta schemes} \label{se a learning method}

{\color{vert}
In this section, we first give a representation of the Runge-Kutta schemes in Definition \ref{de RK scheme} as the solution of a sequence of optimisation problems.  This new representation leads naturally to the use of non-linear regression to compute the scheme in practice. This is an extension of the method proposed in  \cite{hure2020deep} for the implicit Euler scheme. We are able here to work with more efficient schemes, in terms of discrete-time error. We also study the error associated to the approximation space based on Neural Networks. In particular, we show that at each stage of computation the error can be made arbitrary small by relying on the universal approximation theorem, stated in Theorem \ref{universal approximation theorem}. We then present the schemes that will be studied numerically in the next section, namely: the Euler schemes, the Crank-Nicolson scheme, two-stage and three-stage explicit Runge-Kutta schemes. We also give the convergence result for these schemes.
}



\subsection{Implementation of Runge-Kutta schemes}
\label{subse implementation RK schemes}
{\color{vert}
Our approach to implement the Runge-Kutta schemes using non-linear regression is inspired by the approach developed in \cite{hure2020deep} for the implicit Euler scheme (named \textbf{DBDP1} scheme therein). However, we need to introduce extra-terms in the computations as soon as we go away from the Euler scheme (see the $A$-terms below).
Moreover, we  need to compute each stage iteratively.

}

\vspace{2mm}
\noindent The whole procedure is based on the following key observation.
\begin{Lemma} \label{le Lemma optim RK total}
For $1 < q \le Q+1$ and $0\leq n\leq N-1 $, define 
\begin{align}\label{eq definition A_nq RK qN}
	A_{n,q} &= \EFp{t_{n,q}}{  \sum_{k=1}^{q-1}  \left( a_{qk}  H_q^n -\alpha_{qk}  H_{q,k}^n \right) h f(X_{n,k},Y_{n,k}, Z_{n,k})   }  ,
\end{align}
where $( Y_{n, k}, Z_{n, k})$ are defined by \eqref{eq de interm stage Y}-\eqref{eq de interm stage Z}. Then, the transition from step $n$ to $n-1$ in the scheme given in Definition \ref{de RK scheme} is solution of the following  optimisation problem 
\begin{align}
(Y_{n,q},Z_{n,q},A_{n,q})&= \mathrm{argmin}_{(\mfy,\mfz,\mfa) \in \cL^2(\cF_{t_{n,q}})} {L}_{n,q}(\mfy,\mfz,\mfa),
\end{align}
with, for some fixed $\mathfrak{c}>0$,
\begin{align}
 L_{n,q}(\mfy,\mfz,\mfa) & :=  \esp{\Big|Y_{n+1} + h \sum_{k=1}^{q-1} a_{qk} f(X_{n,k},Y_{n,k}, Z_{n,k}) -  \set{\mfy -ha_{qq} f(X_{n,q},\mfy, \mfz)+(\mfz+\mfa)\Delta W_{n, q}  }\Big|^2 
\nonumber \right. \\
&\left. 
+h\mathfrak{c} \Big|\mfa - \sum_{k=1}^{q-1}  \left( a_{qk}  H_q^n -\alpha_{qk}  H_{q,k}^n \right) h f(X_{n,k},Y_{n,k}, Z_{n,k}) \Big|^2}. 
\nonumber
\end{align}

\end{Lemma}

\proof
We first observe that from \eqref{eq de interm stage Y}, \eqref{eq de interm stage Z} and \eqref{eq definition A_nq RK qN} that
\begin{align}\label{eq process Y RK general}
 Y_{n+1} = Y_{n, q} - h \sum_{k=1}^q a_{qk} f(X_{n,k},Y_{n,k}, Z_{n,k}) + (Z_{n, q} + A_{n, q}) \Delta W_{n, q}  - \Delta M_{n, q} ,
 \end{align}
where $\EFp{t_n}{\Delta M_{n, q}}=\EFp{t_n}{\Delta M_{n, q} \Delta W_{n, q}  } = 0, \EFp{t_n}{|\Delta M_{n, q}|^2} < \infty$.  We also have that
\begin{align*}
&\quad\; \esp{\left| \mfa - \sum_{k=1}^{q-1}  \left( a_{qk}  H_q^n -\alpha_{qk}  H_{q,k}^n \right) h f(X_{n,k},Y_{n,k}, Z_{n,k}) \right|^2 } \\
 & = \esp{ | \mfa - A_{n, q} |^2 + \left| A_{n, q} - \sum_{k=1}^{q-1}  \left( a_{qk}  H_q^n -\alpha_{qk}  H_{q,k}^n \right) h f(X_{n,k},Y_{n,k}, Z_{n,k}) \right|^2 }\,.
\end{align*}
Inserting \eqref{eq process Y RK general} into the definition of $L_{n, q}(y,z,a)$ and using the previous equality, we compute
\begin{align*}
&\quad L_{n, q}(\mfy,\mfz, \mfa) \\
 &= \esp{ \Big|Y_{n+1} + h \sum_{k=1}^{q-1} a_{qk} f(X_{n,k},Y_{n,k}, Z_{n,k}) -  \set{\mfy -ha_{qq} f(X_{n,q},\mfy, \mfz)+(\mfz+\mfa)\Delta W_{n, q}  }\Big|^2  } \\
& + \mathfrak{c} h\esp{ | \mfa - A_{n, q} |^2 + \left| A_{n, q} - \sum_{k=1}^{q-1}  \left( a_{qk}  H_q^n -\alpha_{qk}  H_{q,k}^n \right) h f(X_{n,k},Y_{n,k}, Z_{n,k}) \right|^2 }
\\
&=\tilde{L}^{n, q}_1 (\mfy,\mfz) + \tilde{L}^{n, q}_{2}(\mfa, \mfz)  + \tilde{L}^{n, q}_{3}(\mfa) + \ell_{n, q},
\end{align*}
where 
\begin{align}
	&\tilde{L}^{n, q}_{1}(\mfy,\mfz) := \esp{|Y_{n, q}  - a_{qq} h f(X_{n, q},Y_{n, q},Z_{n, q} ) -\set{\mfy -  a_{qq} h f(X_{n, q},\mfy,\mfz)}|^2}, \label{eq loss y}\\
	& \tilde{L}^{n, q}_{2}(\mfa, \mfz) := c_q h \esp{|Z_{n, q} - z + A_{n, q} -\mfa|^2}, \\
	&  \tilde{L}^{n, q}_{3}(\mfa) := \mathfrak{c} h\esp{ | A_{n, q} - \mfa |^2 }, \\
	 & \ell_{n, q} = \esp{|\Delta M_{n, q}|^2 + \mfc h  \left| A_{n, q} - \sum_{k=1}^{q-1}  \left( a_{qk}  H_q^n -\alpha_{qk}  H_{q,k}^n \right) h f(X_{n,k},Y_{n,k}, Z_{n,k}) \right|^2 }.
\end{align}
We then observe that $\tilde{L}^{n, q}_{3}(A_n) =0 $, $\tilde{L}^{n, q}_{2}(A_n, Z_n)=0$ and $\tilde{L}^{n, q}_{1}(Y_n,Z_n)=0$, so that $(Y_{n, q},Z_{n, q}, A_{n, q})$ does achieve the minimum of $L_{n, q}$.\\
Reciprocally, any optimal solution $(\mfy^\star,\mfz^\star, \mfa^\star)$ must satisfy $\tilde{L}^{n, q}_{3}(\mfa^\star)= 0 , \tilde{L}^{n, q}_{2}(\mfa^\star, \mfz^\star) = 0$, which implies $\mfa^\star = A_{n, q}, \mfz^\star=Z_{n, q}$.
Moreover, necessarily one has $\tilde{L}^{n, q}_{1}(\mfy^\star, Z_{n, q}) = 0$, which leads to
\begin{align*}
Y_{n, q}  - a_{qq} h f(X_{n, q},Y_{n, q},Z_{n, q} ) -\set{\mfy^\star -  a_{qq} h f(X_{n, q},\mfy^\star,Z_{n, q})} = 0\;.
\end{align*}
Inserting  \eqref{eq de interm stage Y} into the previous equality, we find
\begin{align*}
\mfy^\star =  \EFp{t_{n,q}}{Y_{n+1} +  h \sum_{k=1}^{q-1} a_{qk} f(X_{n,k},Y_{n,k}, Z_{n,k}) + a_{qq}h f(X_{n,q}, \mfy^\star , Z_{n,q})   }\;.
\end{align*}
By uniqueness of the scheme definition, recall condition \eqref{eq cond well posedness of the scheme}, we get $\mfy^\star = Y_{n, q}$, which concludes the proof.
\eproof

\subsubsection{Scheme definition}

{\color{vert}
We now introduce the scheme that will be implemented in the next section. Our starting point is the representation of the scheme as the solution of a sequence of optimisation problem established in the previous section. From a numerical perspective, this optimisation is solved on a parametric space of approximating functions relying on the Markov property satisfied by the Runge-Kutta schemes. We employ neural network as approximating functions, recall Section \ref{subse NN notations}.
}

\vspace{2mm}
\noindent Let $\Phi = (\Phi_1, \cdots, \Phi_{Q+1}) \in \cC(\R^d,\R)^{Q+1}$ and $\Psi = (\Psi_1, \cdots, \Psi_{Q+1}) \in \cC(\R^d,\R^d)^{Q+1}$, we introduce a generic loss function at each stage of computation $(n,q)$, $1\le n<N$, $1 < q \le Q+1$: 
\begin{align}\label{eq de loss generic}
&\mathrm{L_{n,q}^{\!RK}}[\Phi,\Psi](\theta)  := \esp{\Big|\Phi_1(X_{n+1}) + h \sum_{k=1}^{q-1} a_{qk} f(X_{n,k},\Phi_k(X_{n,k}), \Psi_k(X_{n,k})) 
\right.\\
&\left.
\quad  - \set{\cU(X_{n,q}; \theta) -ha_{qq} f(X_{n,q},\cU(X_{n,q}; \theta), \cV(X_{n,q}; \theta))
+(\cV(X_{n,q}; \theta)+\cA(X_{n,q}; \theta))\Delta W_{n, q}   }\Big|^2 
\nonumber \right. \\
&\left. 
\quad + \mfc h \Big|\cA(X_{n,q}; \theta) - \sum_{k=1}^{q-1}  \left( a_{qk}  H_q^n -\alpha_{qk}  H_{q,k}^n \right) h f(X_{n,k},\Phi_k(X_{n,k}), \Psi_k(X_{n,k})) \Big|^2
}, \nonumber
\end{align}
with $(\cU,\cV,\cA) := \cN_m \in \cS_{d_0, d_1, L, m}^{\rho}(\mathbb{R}^{N_m})$, recall
\eqref{eq de NN function}-\eqref{eq set of networks}. 

{\color{vert}
\noindent In the definition of $\mathrm{L_{n,q}^{\!RK}}$ above, we do not indicate, for the reader's convenience, the dependence upon  the ``balance number''  $\mfc$ (which is fixed). 
}

\vspace{2mm}
\begin{Definition}[Implemented Runge-Kutta scheme]\label{de implemented RK scheme} The numerical solution is computed using the following procedure:
\begin{itemize}
	\item For $n = N$, initialize  $ \hat{\cU}_N = g , \hat{\cV}_N =  \sigma^\top \nabla_x g $,
	\item For $n=N-1, \cdots, 0$: given $ (\hat{\cU}_{n+1}, \hat{\cV}_{n+1})=:(\hat{\cU}_{n,1}, \hat{\cV}_{n,1})  $ do
	\begin{itemize}
	\item for $1 < q \le Q+1$ given   and $(\hat{\cU}_{n,k}, \hat{\cV}_{n,k})$, \textcolor{black}{$1 \leq k < q$}:
	\begin{itemize}
		\item set $(\Phi_k,\Psi_k) := (\hat{\cU}_{n,k}, \hat{\cV}_{n,k})$, $1 \le k < q$, $(\Phi_k,\Psi_k) := 0$, $k \ge q$
		\item Compute a minimizer of the loss function:
		\begin{align}\label{eq optim in scheme}
		\theta_{n,q}^\star \in \mathrm{argmin}_{\theta}\;\mathrm{L_{n,q}^{\!RK}}[\Phi,\Psi](\theta),
		\end{align}
		where $\mathrm{L_{n,q}^{\!RK}}$ defined by  \eqref{eq de loss generic}.		 		
	 	\item set $ (\hat{\cU}_{n,q},\hat{\cV}_{n,q}, \hat{\cA}_{n,q}) :=  \cN_m(\cdot; \theta_{n,q}^\star) \in \cS_{d_0, d_1, L, m}^{\rho}(\mathbb{R}^{N_m})$, recall
\eqref{eq de NN function}-\eqref{eq set of networks}.
	\end{itemize}
	\item Set $ (\hat{\cU}_{n}, \hat{\cV}_{n}):=(\hat{\cU}_{n,Q+1}, \hat{\cV}_{n,Q+1})  $.
	\end{itemize}
\end{itemize}
\end{Definition}


{\color{vert}
\begin{Remark}
\begin{enumerate}[(i)]
\item  Let us note  that we do not have any theoretical guarantees that the minimization problems \eqref{eq optim in scheme} are well-posed on the whole space $ \mathbb{R}^{N_m} $. A sufficient (but far from being necessary) condition is the convexity and coercivity in $\theta$ of the objective loss function \eqref{eq de loss generic}. However, in practical implementation,  the domain is often restricted to a compact subset of $ \mathbb{R}^{N_m} $, so that the minimization problems are always well-posed by continuity of the loss function. In the sequel, our theoretical results are obtained by assuming the well-posedness of the scheme.
\item The value of the ``balance number'' $\mfc$, recall \eqref{eq de loss generic}, may impact slightly the numerical results, see \Cref{balance} in \Cref{se numerics}. 
\item In practice the optimisation problems are solved classically using
a Stochastic Gradient Descent (SGD) procedure which is detailed in Section \ref{se numerics}.
\end{enumerate}
\end{Remark}
}

\subsubsection{Pseudo-consistency of the implemented scheme}
We study here a kind of minimal consistency of the implemented scheme in terms of approximation error made at each step. 

{\color{vert}
Our first result below clarifies the error made between one step of the \textcolor{black}{genuine} scheme given in \eqref{eq de interm stage Y}-\eqref{eq de interm stage Z} and one step of the approximation scheme given above. As expected, this error is controlled by the approximation power of the class of DNN considered.
}

\begin{Lemma} \label{le error one stage generic RK}
Let $0\leq n \leq N-1$ and $1 < q \le Q+1$.  Let $\Phi = (\Phi_1, \cdots, \Phi_{Q+1}) \in \cC(\R^d,\R)^{Q+1}$ and $\Psi = (\Psi_1, \cdots, \Psi_{Q+1}) \in \cC(\R^d,\R^d)^{Q+1}$ and define
\begin{align}
{Y}^{(\Phi,\Psi)}_{n,q} &:= \EFp{t_{n,q}}{
{Y}^{(\Phi,\Psi)}_{n,1} + h \sum_{k=1}^{q} a_{qk} f(X_{n,k},{Y}^{(\Phi,\Psi)}_{n,k} , {Z}^{(\Phi,\Psi)}_{n,k} )   }, \label{eq true NN RK y}
\\
{Z}^{(\Phi,\Psi)}_{n, q} &:= \EFp{t_n}{H^n_q{Y}^{(\Phi,\Psi)}_{n,1} +  h \sum_{k=1}^{q-1} \alpha_{qk} H^n_{q,k}f(X_{n,k},{Y}^{(\Phi,\Psi)}_{n,k} , {Z}^{(\Phi,\Psi)}_{n,k} )  }, \label{eq true NN RK z}
\\
A^{(\Phi,\Psi)}_{n,q} &= \EFp{t_{n,q}}{  \sum_{k=1}^{q-1}  \left( a_{qk}  H_q^n -\alpha_{qk}  H_{q,k}^n \right) h f(X_{n,k},{Y}^{(\Phi,\Psi)}_{n,k}, {Z}^{(\Phi,\Psi)}_{n,k})   }  ,
\label{eq true NN RK a}
\end{align}
with $({Y}^{(\Phi,\Psi)}_{n,k},{Z}^{(\Phi,\Psi)}_{n,k}) := (\Phi_k(X_{n,k}), \Psi_k(X_{n,k}))$, for $1 \le k < q$.
Assume that there exists $\theta^\star \in \R^{N_m}$ such that
\begin{align*}
\theta^\star \in \mathrm{argmin}_{\theta}\;\mathrm{L_{n,q}^{\!RK}}[\Phi,\Psi](\theta),
\end{align*}
then, there exists $C<\infty $ independent of $n$ and $q$ such that, for $h$ small enough, 
\begin{align}\label{eq error one step estim cond}
& \quad \esp{|{Y}^{(\Phi,\Psi)}_{n,q} - \cU_{n,q}(X_{n,q};\theta^\star) |^2+ h|{Z}^{(\Phi,\Psi)}_{n,q} - \cV_{n,q}(X_{n,q}; \theta^\star) |^2
+ h|{A}^{(\Phi,\Psi)}_{n,q} - \cA_{n,q}(X_{n,q};\theta^\star) |^2} \nonumber \\
& \le C \cE_{n,q}(\Phi,\Psi),
\end{align}
where
\begin{align}\label{eq main error NN generic}
\cE_{n,q}(\Phi,\Psi)= \epsilon_{n,q}^{\cN, y}(\Phi,\Psi)+h\epsilon_{n,q}^{\cN, z}(\Phi,\Psi)+h\epsilon_{n,q}^{\cN, a}(\Phi,\Psi),
\end{align}
and 
\begin{align}
	\epsilon_{n,q}^{\cN, y}(\Phi,\Psi)&:=  
	\inf_{\theta^y \in N_m^{d_1=1}} \esp{ |{Y}^{(\Phi,\Psi)}_{n,q}  - \cU_{n,q}(X_{n,q}; \theta^y ) |^2 }, \label{eq NN gen RK error y}\\
	\epsilon_{n,q}^{\cN, a}(\Phi,\Psi) &:=  
	\inf_{\theta^a\in N_m^{d_1=d}} \esp{  |{A}^{(\Phi,\Psi)}_{n,q}  - \cA_{n,q}(X_{n,q}; \theta^a ) |^2}, \label{eq NN gen RK error a} \\
	\epsilon_{n,q}^{\cN, z}(\Phi,\Psi)&:=  
	\inf_{\theta^z \in N_m^{d_1=d}} \esp{  |{Z}^{(\Phi,\Psi)}_{n,q} - \cV_{n,q}(X_{n,q}; \theta^z ) |^2}  \label{eq NN gen RK error z}.
\end{align}
\end{Lemma}

\proof
1. We first observe that \eqref{eq true NN RK y} rewrites
\begin{align}\label{eq process Y RK phi,psi general}
 {Y}^{(\Phi,\Psi)}_{n+1} ={Y}^{(\Phi,\Psi)}_{n,q} &- h \sum_{k=1}^{q} a_{qk} f(X_{n,k},{Y}^{(\Phi,\Psi)}_{n,k} , {Z}^{(\Phi,\Psi)}_{n,k} )+ \left({Z}^{(\Phi,\Psi)}_n+ {A}^{(\Phi,\Psi)}_n\right)\Delta W_{n, q}  + \Delta {M}^{(\Phi,\Psi)}_{n, q} \;,
 \end{align}
 where $\EFp{t_{n, q}}{\Delta M^{(\Phi,\Psi)}_{n, q}}=\EFp{t_{n, q}}{\Delta M^{(\Phi,\Psi)}_{n, q} \Delta W_{n, q} } = 0, \EFp{t_{n, q}}{|\Delta M^{(\Phi,\Psi)}_{n, q}|^2} < \infty$. 
 Following the same computations as in the proof of Lemma \ref{le Lemma optim RK total}, we obtain that
\begin{align}
\mathrm{L_{n, q}^{\!RK}}[\Phi,\Psi](\theta)  = 
\tilde{L}^{n,q}_{1}[\Phi,\Psi](\theta) +
\tilde{L}^{n,q}_{2}[\Phi,\Psi](\theta) +
\tilde{L}^{n,q}_{3}[\Phi,\Psi](\theta)  + \ell_{n,q}[\Phi,\Psi],
\end{align}
with
\begin{align}
	\tilde{L}^{n, q}_{1}[\Phi,\Psi](\theta) &= \mathbb{E} \Big[ |{Y}^{(\Phi,\Psi)}_{n, q}  - a_{qq}h f(X_{n, q},{Y}^{(\Phi,\Psi)}_{n, q} ,{Z}^{(\Phi,\Psi)}_{n, q} )   \\ \nonumber
	 & \qquad - \set{\cU_{n, q}(X_{n, q}; \theta ) - a_{qq} h f(X_{n, q},\cU_{n, q}(X_{n, q}; \theta ),\cV_{n, q}(X_{n, q}; \theta ))}|^2 \Big], \\
	 \tilde{L}^{n, q}_{2}[\Phi,\Psi](\theta)  &= c_q h \esp{|{Z}^{(\Phi,\Psi)}_{n, q}-\cV_{n, q}(X_{n, q}; \theta ) + {A}^{(\Phi,\Psi)}_{n, q} -\cA_{n, q}(X_{n, q}; \theta )|^2}, \\
	  \tilde{L}^{n, q}_{3}[\Phi,\Psi](\theta) &= \mfc h\esp{ | {A}^{(\Phi,\Psi)}_{n, q} - \cA_{n, q}(X_{n, q}; \theta ) |^2 }, \\
	  \ell_{n, q}[\Phi,\Psi] = & \esp{   \mfc h \left|    {A}^{(\Phi,\Psi)}_{n, q} - \sum_{k=1}^{q-1}  \left( a_{qk}  H_q^n -\alpha_{qk}  H_{q,k}^n \right) h f(X_{n,k},{Y}^{(\Phi,\Psi)}_{n,k}, {Z}^{(\Phi,\Psi)}_{n,k}) \right|^2     + |\Delta M^{(\Phi,\Psi)}_{n, q}|^2  \label{eq ln phi psi RK}  }. 
\end{align}
Setting $\tilde{L}_{n, q}[\Phi,\Psi](\theta) := \tilde{L}^{n, q}_{1}[\Phi,\Psi](\theta) +
\tilde{L}^{n, q}_{2}[\Phi,\Psi](\theta) +
\tilde{L}^{n, q}_{3}[\Phi,\Psi](\theta) $, we  then deduce  that 
\begin{align}\label{eq equality of theoretical optim RK}
\mathrm{argmin}_\theta\, \tilde{L}_{n, q}[\Phi,\Psi](\theta)= \mathrm{argmin}_\theta\, \mathrm{L_{n, q}^{\!RK}}[\Phi,\Psi](\theta) .
\end{align}
2. From the Lipschitz continuity of $f$ and recalling \eqref{eq conditions on H 1}-\eqref{eq conditions on H 2}, we directly obtain the following upper bound:
\begin{align}\label{eq upper bound RK}
\tilde{L}_{n, q}[\Phi,\Psi](\theta)   \le &
C\esp{|{Y}^{(\Phi,\Psi)}_{n, q} - \cU_{n, q}(X_{n, q};\theta) |^2  + h|{Z}^{(\Phi,\Psi)}_{n, q} - \cV_{n, q}(X_{n, q};\theta) |^2
 \right. \nonumber  \\  & \qquad \left.  + h|{A}^{(\Phi,\Psi)}_{n, q} - \cA_{n, q}(X_{n, q};\theta) |^2}\,.
\end{align}
We now prove a similar lower bound for $\tilde{L}_{n, q}[\Phi,\Psi](\theta)$. We will use the following inequality
\begin{align} \label{eq useful identity RK}
(x+y)^2 \ge x^2(1-\alpha) + y^2(1- \frac{1}{\alpha}), \quad 0<\alpha<1.
\end{align}  
Thus, for any $\alpha$ such that $1 > \alpha \textcolor{black}{> \frac{2c_q}{2c_q+ \mfc}}$, we obtain
\begin{align}\label{eq interm tilde L RK}
\tilde{L}_2^{n, q}[\Phi,\Psi](\theta) + \tilde{L}_3^{n, q}[\Phi,\Psi](\theta) \ge & h \textcolor{black}{(1-\alpha) c_q}\esp{|{Z}^{(\Phi,\Psi)}_{n, q} -\cV_{n, q}(X_{n, q}; \theta )|^2}
 \\ \nonumber  & 
+h\frac{\textcolor{black}{\mfc}}2\esp{ | {A}^{(\Phi,\Psi)}_{n, q} - \cA_{n, q}(X_{n, q}; \theta ) |^2 }\,.
\end{align}
Using again \eqref{eq useful identity RK} with $\alpha=1/2$ and then the fact that $f$ is Lipschitz continuous, we get
\begin{align*}
\tilde{L}_1^{n, q}[\Phi,\Psi](\theta) &\ge \frac12\esp{|{Y}^{(\Phi,\Psi)}_{n, q} - \cU_{n, q}(X_{n, q},\theta) |^2}
\\
&- a_{qq}^2 h^2 \esp{|f(X_{n, q},{Y}^{(\Phi,\Psi)}_{n, q} ,{Z}^{(\Phi,\Psi)}_{n, q} ) - f(X_{n, q},\cU_{n, q}(X_{n, q}; \theta ),\cV_{n, q}(X_{n, q}; \theta ))|^2}\\
&  \ge 
 (\frac12- 2a_{qq}^2 L^2h^2)\esp{|{Y}^{(\Phi,\Psi)}_{n, q} - \cU_{n, q}(X_{n, q},\theta) |^2} 
-2a_{qq}^2 L^2h^2\esp{|{Z}^{(\Phi,\Psi)}_{n, q} - \cV_{n, q}(X_{n, q},\theta) |^2}.
\end{align*}
Combining the previous inequality with \eqref{eq interm tilde L RK}, we deduce that  for $h$ small enough, 
{
\begin{align}
 & \quad \esp{ |{Y}^{(\Phi,\Psi)}_{n, q} - \cU_{n, q}(X_{n, q}; \theta) |^2+h|{Z}^{(\Phi,\Psi)}_{n, q} - \cV_{n, q}(X_{n, q}; \theta) |^2 +
  h| {A}^{(\Phi,\Psi)}_{n, q} - \cA_{n, q}(X_{n, q}; \theta ) |^2 }  \nonumber \\
  & \le C \tilde{L}_{n, q}[\Phi,\Psi](\theta)
\end{align}
}
3. The above inequality is \emph{a fortiori} true at the optimum $\theta^\star$. Moreover, optimizing on separated networks is always more costly than optimizing on a fully connected network thus leading to  \eqref{eq error one step estim cond}.
\eproof

\vspace{2mm}
\noindent 
{\color{vert} We conclude this section with a result expressing the global error between the scheme given in Definition \ref{de RK scheme} and the one given in Definition \ref{de implemented RK scheme}.
}

\begin{Proposition} \label{pr pseudo consistency gen RK}
Assume that the scheme given in Definition \ref{de implemented RK scheme}  is {well-posed} then
\begin{align}
\max_{0\leq n \le N}\esp{|Y_n - \hat{\cU}_n(X_n)|^2} + h\sum_{n=0}^{N}\esp{|Z_n - \hat{\cV}_n(X_n)|^2}
\le 
C N \sum_{n=0}^{N-1} \bar{\cE}_n\,,
\end{align}
where 
$\bar{\cE}_n := \sum\limits_{k=1}^{Q}\cE_{n,k}(\hat{\Phi}_n,\hat{\Psi}_n)$ recalling \eqref{eq main error NN generic} and with $\hat{\Phi} = (\hat{\cU}_{n,k})_{1 \le k \le Q}$ and $\hat{\Psi} =  (\hat{\cV}_{n,k})_{1 \le k \le Q}$.
\end{Proposition}

\proof  
Let us define, for $0\leq n\leq N-1, 1<q\le Q+1$, recalling \eqref{eq true NN RK y}-\eqref{eq true NN RK z}-\eqref{eq true NN RK a},
\begin{align}
\bar{Y}_{n, q} = Y^{(\hat{\Phi},\hat{\Psi} )}_{n, q} \,,\;
\bar{Z}_{n, q} = Z^{(\hat{\Phi},\hat{\Psi} )}_{n, q} \,,\;
\bar{A}_{n, q} = A^{(\hat{\Phi},\hat{\Psi} )}_{n, q} \,.
\end{align}
We first observe that $(\hat{U}_{n, q},\hat{V}_{n,q}):=(\textcolor{black}{\hat {\cU}_{n,q}(X_{n,q}), \hat{\cV}_{n,q}(X_{n,q})})$ can be rewritten as a perturbed scheme, namely
 \begin{align}
	 	\hat{U}_{n,q} &= \EFp{t_{n,q}}{ \hat{U}_{n+1} + h \sum_{k=1}^{q} a_{qk} f(X_{n,k},\hat{U}_{n,k}, \hat{V}_{n,k}) } + \zeta^y_{n,q},
		\\
	 	\hat{V}_{n,q} &= \EFp{t_{n,q}}{H_q^n \hat{U}_{n+1} + h \sum_{k=1}^{q-1} \alpha_{qk} H_{q,k}^n f(X_{n,k},\hat{U}_{n,k}, \hat{V}_{n,k}) } + \zeta^z_{n,q},
	 \end{align}
with
\begin{align}
\zeta^y_{n,q} &:= \hat{U}_{n,q} - \bar{Y}_{n,q} + a_{qq} h \left( f(X_{n,q},\bar{Y}_{n,q},\bar{Z}_{n,q}) - f(X_{n,q},\hat{U}_{n,q},\hat{V}_{n,q}) \right) , \\
\zeta^z_{n,q} &:= \hat{V}_{n,q} - \bar{Z}_{n,q} \,.
\end{align}
Indeed, with our notations, we have $(\hat{U}_{n, k},\hat{V}_{n, k})= (Y_{n, k}^{(\hat{\Phi},\hat{\Psi})} , Z_{n, k}^{(\hat{\Phi},\hat{\Psi})}), 1\le k < q$.
Moreover, since $\hat{U}_{n,q} = \hat{\cU}_{n,q}(X_{n,q}) = {\cU_{n,q}}(X_{n,q},\theta^\star_{n,q})$, recall Definition \ref{de implemented RK scheme}, it holds
\begin{align*}
 \esp{\frac{1}{h} |\zeta^y_{n, q}|^2+h|\zeta^z_{n, q}|^2}  &\le C\esp{ \frac{1}{h} |{\cU_{n, q}}(X_{n, q},\theta^\star_{n, q}) - \bar{Y}_{n, q}|^2 + h|{\cV_{n, q}}(X_{n, q},\theta^\star_{n, q}) - \bar{Z}_{n, q}|^2}
\\
&\le C N\cE_{n, q}(\hat{\cU}_{n+1},\hat{\cV}_{n+1}),
\end{align*}
where for the last inequality  we applied Lemma \ref{le error one stage generic RK}. Now the proof is concluded using the stability result given in Proposition \ref{pr generic stab for cond exp estim}.
\eproof


{\color{vert}
\vspace{2mm}
\noindent We conclude this section with the following result.
\begin{Theorem} \label{th control global error}Recall the definition of $(\bar{Y}_n,\bar{Z}_n)_{0 \le n \le N}$ in \eqref{eq de support bsde} and of the discrete-time error $\cT_N$ in \eqref{eq de truncation error}. Then, the following holds
\begin{align}\label{eq control global}
\max_{n \le N}\esp{|\bar{Y}_n - \hat{\cU}_n(X_n)|^2} + h\sum_{n=0}^{N}\esp{|\bar{Z}_n - \hat{\cV}_n(X_n)|^2}
\le
C(\cT_N +  N \sum_{n=0}^{N-1} \bar{\cE}_n)\,,
\end{align}
where 
$\bar{\cE}_n := \sum\limits_{k=1}^{Q}\cE_{n,k}(\hat{\Phi}_n,\hat{\Psi}_n)$ and  with $\hat{\Phi} = (\hat{\cU}_{n,k})_{1 \le k \le Q}$ and $\hat{\Psi} =  (\hat{\cV}_{n,k})_{1 \le k \le Q}$, recall \eqref{eq main error NN generic}.
\end{Theorem}
\proof We simply observe that
\begin{align*}
&\quad  \max_n \esp{|\bar{Y}_n-\hat{\cU}_n(X_n)|^2} + \sum_{n=0}^{N-1}h\esp{|\bar{Z}_n-\hat{\cV}_n(X_n)|^2}\\
& \le 2\left(\max_n \esp{|\bar{Y}_n-Y_n|^2} + \sum_{n=0}^{N-1}h\esp{|\bar{Z}_n-Z_n|^2}\right)
\\
&+ 2\left(\max_n \esp{|{Y}_n-\hat{\cU}_n(X_n)|^2} + \sum_{n=0}^{N-1}h\esp{|{Z}_n-\hat{\cV}_n(X_n)|^2}\right).
\end{align*}
The first term in the right hand side of the previous inequality is the discrete-time error and the second term is upper bounded using Proposition \ref{pr pseudo consistency gen RK}.
\eproof
}

{\color{vert}
\begin{Remark}
\begin{enumerate}[(i)]
\item Observe that \eqref{eq control global} provides an upper-bound for the error $|\cY_0 - \hat{\cU}_n(X_0)|$ at the initial time.
\item The above result is a first step only toward the proof of convergence of the implemented RK scheme. Indeed, one observes that the local error $\bar{\cE}_n$ is defined using the solution computed at the previous time step. For fixed number of time steps $N$, Theorem \ref{universal approximation theorem} combined with Lemma \ref{le error one stage generic RK} shows that this error can be made as small as desired. This confirms that the implemented scheme is a reasonable scheme to consider. The theoretical difficulty to obtain a full convergence result comes from the highly non linear optimisation problem that have to be solved when using deep neural networks as approximation space. From a practical point of view, next section shows that the numerical procedure associated to Definition \ref{de implemented RK scheme} is efficient.
\end{enumerate}
\end{Remark}
}

\subsection{Review of classical Runge-Kutta schemes}
\label{subse implemented schemes}
{\color{vert} We now present  the schemes that are used to perform the numerical simulations. We will consider only the Euler schemes, the Crank-Nicholson scheme and some two-stage and three-stage explicit Runge-Kutta schemes, recall Remark \ref{re restriction}.
}

\subsubsection{Euler schemes}
\label{se euler scheme basics}
{\color{vert}
We first specialize the results of the previous section to the case of the Euler schemes. As already mentioned, for the implicit Euler scheme, we obtain exactly the procedure introduced in \cite{hure2020deep}.
}


\noindent For this section, we assume that $(X_n)_{n\le N}$ is given by the classical Euler scheme on $\pi$:
\begin{align}
X_{n+1} = X_n + \mu(X_n)h + \sigma(X_n)\Delta W_n\;, n \le N-1\;,
\end{align}
with $\Delta W_n := W_{t_{n+1}} - W_{t_n} $ and $X_0 = \cX_0$.

\noindent The implicit Euler scheme for BSDEs \cite{zhang2004numerical, bouchard2004discrete}, which  is a one stage scheme, reads classically as follows, 
for $0 \le n <N$,
\begin{align}\label{imp Euler scheme BSDEs func}
Y_n= \EFp{t_n}{Y_{n+1}+ h f(X_n,Y_n,Z_n) )} \text{ and } Z_n = \EFp{t_n}{\frac{\Delta W_n}{h} Y_{n+1} },
\end{align}
where we set $H_n = \frac{W_{t_{n+1}} - W_{t_n}}{h}$. 
\textcolor{vert}{We observe that  \eqref{eq conditions on H 1}-\eqref{eq conditions on H 2} and \eqref{eq conditions on H 3} are satisfied and  that $\HYP{X}_1$ holds true here.} 

\vspace{2mm}
\noindent The loss function given in \eqref{eq de loss generic} simplifies here as follows: for $\varphi \in \cC(\R^d,\R)$,
\begin{align}\label{eq loss euler imp}
\mathrm{L_n^{\!EUi}}[\varphi](\theta) & := \esp{ \Big| \varphi(X_{n+1}) - \big\{ \cU(X_{n}; \theta)  - h f(X_n, \cU(X_{n}; \theta), {\cV}(X_{n}; \theta) )    +  {\cV}(X_{n}; \theta)  \Delta W_n  \big\}   \Big|^2  },
\end{align}
with $(\cU,\cV) := \cN_m \in \cS_{d_0, d_1, L, m}^{\rho}(\mathbb{R}^{N_m})$, recall
\eqref{eq de NN function}-\eqref{eq set of networks}.

\noindent The implemented implicit Euler scheme is then given by

\begin{Definition}[Implemented implicit Euler scheme]\label{de euler imp} The numerical solution is computed using the following step:
\begin{itemize}
	\item For $n = N$, initialize  $ \hat{\cU}_N = g , \hat{\cV}_N =  \sigma^\top \nabla_x g $. 
	\item For $n=N-1, \cdots, 1, 0$, given $ \hat{\cU}_{n+1} $, 
	\begin{itemize}
		\item Compute a minimizer of the loss function:
		\begin{align}\label{opti euler imp defi}
		\theta_n^\star \in \mathrm{argmin}_{\theta}\;\mathrm{L_n^{\!EUi}}[\hat{\cU}_{n+1}](\theta),
		\end{align}
		recall \eqref{eq loss euler imp}.		 		
	 	\item Set $ (\hat{\cU}_{n},\hat{\cV}_{n}) :=  \cN_m(\cdot; \theta_n^\star) \in \cS_{d_0, d_1, L, m}^{\rho}(\mathbb{R}^{N_m})$, recall
\eqref{eq de NN function}-\eqref{eq set of networks}.
	\end{itemize}
\end{itemize}
\end{Definition}

\vspace{2mm}

We will also consider 
the explicit Euler scheme which has essentially the same empirical results, see \Cref{se numerics}. The theoretical scheme writes
\begin{align}\label{exp Euler scheme BSDEs func}
Y_n= \EFp{t_n}{Y_{n+1}+ h f(X_n,Y_{n+1},Z_{n+1}) )} \text{ and } Z_n = \EFp{t_n}{\frac{\Delta W_n}{h} Y_{n+1} },
\end{align}
for $0\le n <N$. The loss function at step $n<N$ is given by: for $(\varphi,\psi) \in \cC(\R^d,\R)\times \cC(\R^d,\R^d)$,
\begin{align}\label{eq loss euler exp}
\mathrm{L_n^{\!EUe}}[\varphi,\psi](\theta) & := \mathbb{E} \Big[ \Big| \varphi(X_{n+1}) + h f(X_{n+1}, \varphi(X_{n+1}), \psi(X_{n+1}))   - \big\{ \cU(X_{n}; \theta)    +  {\cV}(X_{n}; \theta)  \Delta W_n  \big\}   \Big|^2  \Big]
\end{align}
with $(\cU,\cV) := \cN_m \in \cS_{d_0, d_1, L, m}^{\rho}(\mathbb{R}^{N_m})$, recall
\eqref{eq de NN function}-\eqref{eq set of networks}.


\begin{Definition}[Implemented explicit Euler scheme]\label{de euler exp imp} The numerical solution is computed using the following step:
\begin{itemize}
	\item For $n = N$, initialize  $ \hat{\cU}_N = g , \hat{\cV}_N =  \sigma^\top \nabla_x g $, $\hat{\cA}_N=0$. 
	\item For $n=N-1, \cdots, 1, 0$, given $ \hat{\cU}_{n+1}, \hat{\cV}_{n+1}  $, 
	\begin{itemize}
		\item Compute a minimizer of the loss function:
		\begin{align}\label{opti euler exp defi}
		\theta_n^\star \in \mathrm{argmin}_{\theta}\;\mathrm{L_n^{\!EUe}}[\hat{\cU}_{n+1},\hat{\cV}_{n+1}](\theta),
		\end{align}
		recall \eqref{eq loss euler exp}.		 		
	 	\item Set $ (\hat{\cU}_{n},\hat{\cV}_{n}) :=  \cN_m(\cdot; \theta_n^\star) \in \cS_{d_0, d_1, L, m}^{\rho}(\mathbb{R}^{N_m})$, recall
\eqref{eq de NN function}-\eqref{eq set of networks}.
	\end{itemize}
\end{itemize}
\end{Definition}

{\color{vert}
\noindent Combining Theorem \ref{th disc error control} and \ref{th control global error}, we can state the following result for the Euler schemes.
\begin{Corollary}
Let $(\bar{Y}_n, \bar{Z}_n)_{n \le N}$ be given by \eqref{eq de support bsde} and assume that $(\mathrm{H}r)_1$ is in force.
Then, the following holds, for $h$ small enough,
\begin{align*}
\max_n \esp{|\bar{Y}_n-\hat{\cU}_n(X_n)|^2} + \sum_{n=0}^{N-1}h\esp{|\bar{Z}_n-\hat{\cV}_n(X_n)|^2}
\le
C(h^{2} + \sum_{n=0}^{N-1} \bar{\cE}^{\mathrm{EU}}_n)\;.
\end{align*}
with $\bar{\cE}^{\mathrm{EU}}_n := \cE_{n,1}\left(\hat{\cU}_{n+1}\right)$ for the implicit Euler scheme and $\bar{\cE}^{\mathrm{EU}}_n := \cE_{n,1}\left((\hat{\cU}_{n+1},\hat{\cV}_{n+1})\right)$ for the explicit Euler scheme, recall \eqref{eq main error NN generic}. 
\end{Corollary}
}

\subsubsection{Crank-Nicolson scheme} \label{subse CN scheme}
We now consider the Crank-Nicolson scheme for BSDEs, which is a one stage scheme. It has been introduced in \cite{crisan2014second}, where it is implemented  using cubature methods and tree based branching algorithm(TBBA).

The theoretical scheme reads as follows. For the $Y-$part, it is the usual Crank-Nicolson scheme, namely 
\begin{align} \label{eq approx Y total}
 \left\{
	 		\begin{array}{rl}
Y_N & = g(X_N),  \\
Y_n & = \EFp{t_n}{Y_{n+1} + \frac{h}2 (f(X_n,Y_n,Z_n)+ f(X_{n+1},Y_{n+1},Z_{n+1}) )}, \quad 0\le n \le N-1,
	\end{array}
	 \right.
\end{align}
and  for the $Z-$part, 
\begin{align} \label{eq approx Z total}
	\left\{
			\begin{array}{l}
				Z_N = \sigma^\top \nabla_x g(X_N), \\
				Z_n = \EFp{t_n}{ H_n(Y_{n+1} + h f(X_{n+1},Y_{n+1},Z_{n+1}) ) },
			\end{array}
	\right.
\end{align} 
where $ H_n \in \mathbb{R}^d  $ is a  $\cF_{t_{n+1}}$-mesurable  random variable, satisfying \eqref{eq conditions on H 1}-\eqref{eq conditions on H 2}. 

\vspace{2mm}
\noindent 
The loss function given in \eqref{eq de loss generic} rewrites as follows in this context. Let $(\varphi,\psi)\in \cC(\R^d,\R)\times \cC(\R^d,\R^d)$, for step $n<N$, we define
\begin{align}\label{eq de loss}
\mathrm{L_n^{\!CN}}[\varphi,\psi](\theta) & := \mathbb{E} \Big[ \Big| \varphi(X_{n+1}) - \big\{ \cU(X_{n}; \theta)   -\frac{h}{2} f( X_{n+1},\varphi(X_{n+1}), \psi(X_{n+1}))  \\ 
	 						& - \frac{h}{2} f(X_n, \cU(X_{n}; \theta), {\cV}(X_{n}; \theta) )    + ({\cV}(X_{n}; \theta) + \cA(X_n; \theta ) ) \Delta W_n \big\}   \Big|^2  \nonumber \\
	 						& + C_0 h\Big|  \frac{h}{2}  f(X_{n+1},\varphi(X_{n+1}), \psi(X_{n+1}) )     H_n  +  \cA (X_n; \theta)  \Big|^2 \Big],
							 \nonumber
\end{align}
with $(\cU,\cV,\cA) := \cN_m \in \cS_{d_0, d_1, L, m}^{\rho}(\mathbb{R}^{N_m})$, recall
\eqref{eq de NN function}-\eqref{eq set of networks}.



\noindent We then consider:
\begin{Definition}[Implemented Crank-Nicolson scheme]\label{de crank-nicolson imp} The numerical solution is computed using the following step:
\begin{itemize}
	\item For $n = N$, initialize  $ \hat{\cU}_N = g , \hat{\cV}_N =  \sigma^\top \nabla_X g $, $\hat{\cA}_N=0$. 
	\item For $n=N-1, \cdots, 1, 0$, given $ \hat{\cU}_{n+1}, \hat{\cV}_{n+1}  $, 
	\begin{itemize}
		\item Compute a minimizer of the loss function:
		\begin{align} \label{opti CN defi}
		\theta_n^\star \in \mathrm{argmin}_{\theta}\;\mathrm{L_n^{\!CN}}[\hat{\cU}_{n+1},\hat{\cV}_{n+1}](\theta),
		\end{align}
		recall \eqref{eq de loss}.		 		
	 	\item Set $ (\hat{\cU}_{n},\hat{\cV}_{n}, \hat{\cA}_n) :=  \cN_m(\cdot; \theta_n^\star) \in \cS_{d_0, d_1, L, m}^{\rho}(\mathbb{R}^{N_m})$, recall
\eqref{eq de NN function}-\eqref{eq set of networks}.
	\end{itemize}
\end{itemize}
\end{Definition}

\begin{Remark}\label{remark of the optim lemma CN total}
\begin{enumerate}[(i)]
\item In the numerics, we  use the following loss function instead of \eqref{eq de loss} in order to reduce the variance of $A-$part by a control variate technique. 
\begin{align}
\mathrm{L_{n, R}^{\!CN}}[\varphi,\psi](\theta) & := \mathbb{E} \Big[ \Big| \varphi(X_{n+1}) - \big\{ \cU(X_{n}; \theta)   -\frac{h}{2} f( X_{n+1},\varphi(X_{n+1}), \psi(X_{n+1}))  \\ 
	 						& - \frac{h}{2} f(X_n, \cU(X_{n}; \theta), {\cV}(X_{n}; \theta) )    + ({\cV}(X_{n}; \theta) + \cA(X_n; \theta ) ) \Delta W_n \big\}   \Big|^2  \nonumber \\
	 						& + \mfc h\Big|  \frac{h}{2} \left(  f(X_{n+1},\varphi(X_{n+1}), \psi(X_{n+1}) ) -  f(X_{n},\varphi(X_{n}), \psi(X_{n}) ) \right)  H_n   + \cA (X_n; \theta)  \Big|^2 \Big].
							 \nonumber
\end{align} 

\end{enumerate}
\end{Remark}

\vspace{2mm}
{\color{vert} 
\noindent We conclude this section with a control of the global error for the Crank-Nicolson scheme combining Theorem \ref{th disc error control} and \ref{th control global error}.
\begin{Corollary}
Let $(\bar{Y}_n, \bar{Z}_n)_{n \le N}$ be given by \eqref{eq de support bsde} and assume that $(\mathrm{H}r)_2$ and $(\mathrm{H}X)_2$ are in force.
Then, the following holds, for $h$ small enough,
\begin{align*}
\max_n \esp{|\bar{Y}_n-\hat{\cU}_n(X_n)|^2} + \sum_{n=0}^{N-1}h\esp{|\bar{Z}_n-\hat{\cV}_n(X_n)|^2}
\le
C(h^{4} + \sum_{n=0}^{N-1} \bar{\cE}^{\mathrm{CN}}_n)\;.
\end{align*}
with $\bar{\cE}^{\mathrm{CN}}_n := \cE_{n,1}\left((\hat{\cU}_{n+1},\hat{\cV}_{n+1})\right)$, recall \eqref{eq main error NN generic}. 
\end{Corollary}
}

\subsubsection{Two stage explicit Runge-Kutta scheme}\label{subse two stage exp RK}
We now present the numerical procedure to compute two stage Runge-Kutta scheme. 
The first stage  is an explicit Euler step and then there is no need to introduce the correction term $A$. 
Recalling Theorem \ref{th disc error control}(iii), one can choose the coefficients such that 
\[ a_{21} = \alpha_{21} = c_2, \quad a_{31} = \alpha_{31} = 1 - \frac{1}{2c_2}, \quad a_{32} = \alpha_{32} = \frac{1}{2c_2},  \]
to obtain the optimal bound on the discrete time error.

\noindent The scheme reads thus as follows
\begin{align}
	 	Y_{n, 2} &= \EFp{t_{n,2}}{ Y_{n+1} +  c_2 h f(X_{n+1},Y_{n+1}, Z_{n+1}) },\\
	 	Z_{n, 2} &= \EFp{t_{n,2}}{H^n_2Y_{n+1} +  H^n_{2,1}c_2 h f(X_{n+1},Y_{n+1}, Z_{n+1}) ) },
\end{align}
and 
\begin{align}
	 	Y_{n} &= \EFp{t_{i}}{ Y_{n+1} +  (1 - \frac{1}{2c_2}) h f(X_{n+1},Y_{n+1}, Z_{n+1}) +  \frac{1}{2c_2}h f(X_{n,2},Y_{n,2}, Z_{n, 2})},\\
	 	Z_{n} &= \EFp{t_{i}}{H^n_3 Y_{n+1} +  (1 - \frac{1}{2c_2}) h H^n_{3} f(X_{n+1},Y_{n+1}, Z_{n+1}) + 
		\frac{1}{2c_2}h H^n_{3,2}  f(X_{n,2},Y_{n,2}, Z_{n, 2}) }.
\end{align}
Note that we have used $H^n_{3,1}=H^n_3$, which simplifies slightly the  term $A_{n,3}$ below.

\noindent We must consider loss functions for each stage of computations, namely: \\
- First stage: For $(\varphi,\psi)\in \cC(\R^d,\R)\times \cC(\R^d,\R^d)$,
\begin{align*}
\mathrm{L_{n,2}^{\!RK2}}[\varphi,\psi](\theta)  &:= 
\esp{\Big|\varphi(X_{n+1}) + h c_2 f(X_{n+1},\varphi(X_{n+1}), \psi(X_{n+1})) 
\right.\\
&\left.
\qquad  
- \left\{\cU(X_{n,2}; \theta) 
+ \cV(X_{n,2}; \theta) \Delta W_{n, 2}   \right\}\Big|^2 
} 
\end{align*}
with $(\cU,\cV)  \in \cS_{d_0, d_1, L, m}^{\rho}(\mathbb{R}^{N_m})$, recall
\eqref{eq de NN function}-\eqref{eq set of networks};\\ 
- last stage: For $(\Phi,\Psi) \in \cC(\R^d,\R)^2 \times\cC(\R^d,\R^d)^2 $,
\begin{align}\label{eq de loss RK2 - last stage}
&\mathrm{L_{n,3}^{\!RK2}}[\Phi,\Psi](\theta)  := \\
&\esp{\Big|\Phi_1(X_{n+1}) + h  (1 - \frac{1}{2c_2})  f(X_{n+1},\Phi_1(X_{n+1}), \Psi_1(X_{n+1})) 
 \nonumber
\right.\\
&\left.
\quad + \frac{h}{2c_2} f(X_{n,2},\Phi_2(X_{n,2}), \Psi_2(X_{n,2}))   - \left\{\cU(X_{n}; \theta) 
+ (\cV(X_{n}; \theta)+\cA(X_{n}; \theta) )\Delta W_{n}  \right\}\Big|^2 
\nonumber \right. \\
&\left. 
+\mfc h \Big|\cA(X_{n}; \theta) -  \frac{1}{2c_2}\left( H_3^n -  H_{3,2}^n \right) h f(X_{n,2},\Phi_2(X_{n,2}), \Psi_2(X_{n,2})) \Big|^2
} \nonumber
\end{align}
with $(\cU,\cV,\cA) := \cN_m \in \cS_{d_0, d_1, L, m}^{\rho}(\mathbb{R}^{N_m})$, recall
\eqref{eq de NN function}-\eqref{eq set of networks}.

\vspace{2mm}
\noindent The implemented scheme in then given by
\begin{Definition}
For a given fixed balance number $\mfc>0$ , the algorithm is designed as follows:
\begin{itemize}
	\item For $i = N$, initialize  $ \hat{\cU}_N = g , \hat{\cV}_N = \sigma^\top \nabla_x g $. 
	\item For $i=N-1, \cdots,  0$, given $ \hat{\cU}_{n+1}, \hat{\cV}_{n+1}  $, 
		\begin{itemize}
		\item Compute a minimizer of the loss function at step $(n,2)$: 
		\begin{align*}
		\theta_{n,2}^\star \in 
		\mathrm{argmin}_{\theta}\;
		\mathrm{L_{n,2}^{\!RK2}}[\hat{\cU}_{n+1},\hat{\cV}_{n+1}](\theta)\;.
		\end{align*}
		Set $(\hat{\cU}_{n,2},\hat{\cV}_{n,2}) = (\cU(\cdot,\theta_{n,2}^\star),\cV(\cdot,\theta_{n,2}^\star))$.
	        \item Compute a minimizer of the loss at step $(n,3)$:
	        \begin{align*}
		\theta_{n,3}^\star \in 
		\mathrm{argmin}_{\theta}\;
		\mathrm{L_{n,3}^{\!RK2}}[(\hat{\cU}_{n+1},\hat{\cU}_{n,2}),(\hat{\cV}_{n+1},\hat{\cV}_{n,2})](\theta)\;.
		\end{align*}
		Set $(\hat{\cU}_{n},\hat{\cV}_{n}) = (\cU(\cdot,\theta_{n,3}^\star),\cV(\cdot,\theta_{n,3}^\star))$.
	\end{itemize}
\end{itemize}
\end{Definition}

{\color{vert}
\noindent Our convergence result for this scheme, obtained by combining Theorem \ref{th disc error control} and Theorem \ref{th control global error}, reads as follows.
\begin{Corollary}
Let $(\bar{Y}_n, \bar{Z}_n)_{n \le N}$ be given by \eqref{eq de support bsde}.
 Assume $(\mathrm{H}r)_2$ and $(\mathrm{H}X)_2$.  Then, the following holds
\begin{align*}
\max_n \esp{|\bar{Y}_n-\hat{\cU}_n(X_n)|^2} + \sum_{n=0}^{N-1}h\esp{|\bar{Z}_n-\hat{\cV}_n(X_n)|^2}
\le
C(h^{4} + \sum_{n=0}^{N-1} \bar{\cE}^{\mathrm{RK2}}_n)\,,
\end{align*}
with $\bar{\cE}^{\mathrm{RK2}}_n = \cE_{n,2}\left((\hat{\cU}_{n+1},\hat{\cU}_{n,2}),(\hat{\cV}_{n+1},\hat{\cV}_{n,2})\right)$, where $\cE_{n,2}$ is defined by \eqref{eq main error NN generic}. 
\end{Corollary}
}

\subsubsection{Three stage explicit Runge-Kutta scheme}\label{subse three stage exp RK}
The numerical procedure of three stage Runge-Kutta scheme consists simply in one more iteration than the two stage Runge-Kutta scheme. 
Recalling Theorem \ref{th disc error control}(iv), one can choose the coefficients such that 
\begin{align*}
	&  a_{21} = \alpha_{21} = c_2,  \quad a_{31} = \alpha_{31} = \frac{c_3(3c_2 - 3c_2^2 - c_3)}{c_2(2-3c_2)}, \quad
	a_{32} = \alpha_{32} = \frac{c_3(c_3 - c_2)}{c_2(2-3c_2)}, \\
	& a_{41} =  \alpha_{41} = \frac{-3c_3 + 6c_2c_3 -3c_2}{6c_2c_3} ,\quad
	a_{42} =  \alpha_{42} = \frac{3c_3 - 2}{6c_2(c_3 - c_2)}, \quad
	a_{43} =  \alpha_{43} = \frac{2-3c_2}{6c_3(c_3 - c_2)}.
\end{align*}
to obtain the optimal bound on the discrete time error.
\noindent The scheme reads thus as follows
\begin{align}
	 	Y_{n, 2} &= \EFp{t_{n,2}}{ Y_{n+1} +  c_2 h f(X_{n+1},Y_{n+1}, Z_{n+1}) },\\
	 	Z_{n, 2} &= \EFp{t_{n,2}}{H^n_2Y_{n+1} +  H^n_{2,1}c_2 h f(X_{n+1},Y_{n+1}, Z_{n+1}) ) },
\end{align}
\begin{align}
	 	Y_{n, 3} &= \EFp{t_{n,3}}{ Y_{n+1} +  a_{31} h f(X_{n+1},Y_{n+1}, Z_{n+1}) + a_{32} h f(X_{n, 2},Y_{n, 2}, Z_{n, 2}) },\\
	 	Z_{n, 3} &= \EFp{t_{n,3}}{H^n_3 Y_{n+1} +  H^n_{3,1} \alpha_{31} h f(X_{n+1},Y_{n+1}, Z_{n+1}) + H^n_{3,2} \alpha_{32} h f(X_{n,2},Y_{n,2}, Z_{n, 2}) )  },
\end{align}
and 
\begin{align}
	 	Y_{n} &= \EFp{t_{n,4}}{ Y_{n+1} +  a_{41} h f(X_{n+1},Y_{n+1}, Z_{n+1}) + a_{42} h f(X_{n, 2},Y_{n, 2}, Z_{n, 2})  
	 	  + a_{43} h f(X_{n, 3},Y_{n, 3}, Z_{n, 3}) },\\
	 	Z_{n} &= \mathbb{E}_{t_{n,4}}[H^n_4 Y_{n+1} +  H^n_{4,1} \alpha_{41} h f(X_{n+1},Y_{n+1}, Z_{n+1})  \nonumber \\
	 	 &\qquad\quad + H^n_{4,2} \alpha_{42} h f(X_{n,2},Y_{n,2}, Z_{n, 2}) )  + H^n_{4,3} \alpha_{43} h f(X_{n,3},Y_{n,3}, Z_{n, 3}) ) ],
\end{align}
Note that we have used $H^n_{q,1}=H^n_q$ for $ q = 3, 4$, which simplifies slightly the  term $A_{n,3}, A_{n,4}$ below.\\
\noindent We must consider loss functions for each stage of computations, namely: \\
- First stage: For $(\Phi,\Psi) \in \cC(\R^d,\R)\times \cC(\R^d,\R^d) $,
\begin{align*}
\mathrm{L_{n,2}^{\!RK3}}[\Phi,\Psi](\theta)  &:= 
\esp{\Big|\varphi(X_{n+1}) + h c_2 f(X_{n+1},\Phi(X_{n+1}), \Psi(X_{n+1})) 
\right.\\
&\left.
\quad  
- \left\{\cU(X_{n,2}; \theta) 
+ \cV(X_{n,2}; \theta) \Delta W_{n, 2}  \right\}\Big|^2 
} 
\end{align*}
with $(\cU,\cV)  \in \cS_{d_0, d_1, L, m}^{\rho}(\mathbb{R}^{N_m})$, recall
\eqref{eq de NN function}-\eqref{eq set of networks};\\ 
- Second stage: For $(\Phi,\Psi) \in \cC(\R^d,\R)^2 \times\cC(\R^d,\R^d)^2 $,
\begin{align}\label{eq de loss RK3 - second stage}
&\mathrm{L_{n,3}^{\!RK3}}[\Phi,\Psi](\theta)  := \esp{\Big|\Phi_1(X_{n+1}) + h  a_{31} f(X_{n+1},\Phi_1(X_{n+1}), \Psi_1(X_{n+1})) 
\right.\\
&\left.
\quad + a_{32} h f(X_{n,2},\Phi_2(X_{n,2}), \Psi_2(X_{n,2}))   - \left\{\cU(X_{n,3}; \theta) 
+ (\cV(X_{n,3}; \theta)+\cA(X_{n,3}; \theta) )\Delta W_{n, 3}  \right\}\Big|^2 
\nonumber \right. \\
&\left. 
\quad +\mfc h \Big|\cA(X_{n,3}; \theta) -\left(a_{32} H_3^n - \alpha_{32} H_{3,2}^n \right) h f(X_{n,2},\Phi_2(X_{n,2}), \Psi_2(X_{n,2})) \Big|^2
} \nonumber
\end{align}
with $(\cU,\cV,\cA) := \cN_m \in \cS_{d_0, d_1, L, m}^{\rho}(\mathbb{R}^{N_m})$, recall
\eqref{eq de NN function}-\eqref{eq set of networks}.\\
- Third stage: For $(\Phi,\Psi) \in \cC(\R^d,\R)^3 \times\cC(\R^d,\R^d)^3 $,
\begin{align}\label{eq de loss RK3 - second stage}
\mathrm{L_{n,3}^{\!RK3}}[\Phi,\Psi](\theta) & := \esp{\Big|\Phi_1(X_{n+1}) + h  a_{41} f(X_{n+1},\Phi_1(X_{n+1}), \Psi_1(X_{n+1})) 
\right.\\
&\left.
 + a_{42} h f(X_{n,2},\Phi_2(X_{n,2}), \Psi_2(X_{n,2})) + a_{43} h f(X_{n,3},\Phi_3(X_{n,3}), \Psi_3(X_{n,3}))  
\nonumber \right. \\
&\left.  - \left\{\cU(X_{n}; \theta) 
+ (\cV(X_{n}; \theta)+\cA(X_{n}; \theta) )\Delta W_{n}   \right\}\Big|^2 
\nonumber \right. \\
&\left. 
+ \mfc h \Big|\cA(X_{n}; \theta) -\left(a_{42} H_4^n - \alpha_{42} H_{4,2}^n \right) h f(X_{n,2},\Phi_2(X_{n,2}), \Psi_2(X_{n,2}))  
\nonumber \right. \\
&\left.  
-\left(a_{43} H_4^n - \alpha_{43} H_{4,3}^n \right) h f(X_{n,3},\Phi_3(X_{n,3}), \Psi_3(X_{n,3})) \Big|^2
} \nonumber
\end{align}
with $(\cU,\cV,\cA) := \cN_m \in \cS_{d_0, d_1, L, m}^{\rho}(\mathbb{R}^{N_m})$, recall
\eqref{eq de NN function}-\eqref{eq set of networks}.

\vspace{2mm}
\noindent The implemented scheme is then:
\begin{Definition}
For a given fixed balance number $\mfc > 0 $, the algorithm is designed as follows:
\begin{itemize}
	\item For $i = N$, initialize  $ \hat{\cU}_N = g , \hat{\cV}_N = \sigma^\top \nabla_x g $. 
	\item For $i=N-1, \cdots, 1, 0$, given $ \hat{\cU}_{n+1}, \hat{\cV}_{n+1}  $, 
		\begin{itemize}
		\item Compute a minimizer of the loss function at step $(n,2)$: 
		\begin{align*}
		\theta_{n,2}^\star \in 
		\mathrm{argmin}_{\theta}\;
		\mathrm{L_{n,2}^{\!RK3}}[\hat{\cU}_{n+1},\hat{\cV}_{n+1}](\theta)\;.
		\end{align*}
		Set $(\hat{\cU}_{n,2},\hat{\cV}_{n,2}) = (\cU(\cdot,\theta_{n,2}^\star),\cV(\cdot,\theta_{n,2}^\star))$.
	        \item Compute a minimizer of the loss at step $(n,3)$:
	        \begin{align*}
		\theta_{n,3}^\star \in 
		\mathrm{argmin}_{\theta}\;
		\mathrm{L_{n,3}^{\!RK3}}[(\hat{\cU}_{n+1},\hat{\cU}_{n,2}),(\hat{\cV}_{n+1},\hat{\cV}_{n,2})](\theta)\;.
		\end{align*}
		Set $(\hat{\cU}_{n},\hat{\cV}_{n}) = (\cU(\cdot,\theta_{n,3}^\star),\cV(\cdot,\theta_{n,3}^\star))$.
		\item Compute a minimizer of the loss at step $(n,4)$:
	        \begin{align*}
		\theta_{n,4}^\star \in 
		\mathrm{argmin}_{\theta}\;
		\mathrm{L_{n,4}^{\!RK3}}[(\hat{\cU}_{n+1},\hat{\cU}_{n,2},\hat{\cU}_{n,3}),(\hat{\cV}_{n+1},\hat{\cV}_{n,2}, \hat{\cV}_{n,3})](\theta)\;.
		\end{align*}
		Set $(\hat{\cU}_{n},\hat{\cV}_{n}) = (\cU(\cdot,\theta_{n,4}^\star),\cV(\cdot,\theta_{n,4}^\star))$.
	\end{itemize}
\end{itemize}
\end{Definition}
The convergence result for the above scheme is stated below in Theorem \ref{thm RK general}.

\vspace{1\baselineskip}

{\color{vert}
\noindent We conclude this section by the following convergence result for three stage explicit Runge-Kutta scheme, obtained by combining Theorem \ref{th disc error control} and Theorem \ref{th control global error}.


\begin{Corollary}\label{thm RK general}
Let $(\bar{Y}_n, \bar{Z}_n)_{n \le N}$ be given by \eqref{eq de support bsde}
Assume $(\mathrm{H}r)_3$ and $(\mathrm{H}X)_3$. Then, the following holds
\begin{align*}
\max_n \esp{|\bar{Y}_n-\hat{\cU}_n(X_n)|^2} + \sum_{n=0}^{N-1}h\esp{|\bar{Z}_n-\hat{\cV}_n(X_n)|^2}
\le
C(h^{6} + \sum_{n=0}^{N-1} \bar{\cE}^{\mathrm{RK3}}_n)\;,
\end{align*}
with $\bar{\cE}^{\mathrm{RK3}}_n = \cE_{n,3}\left((\hat{\cU}_{n+1},\hat{\cU}_{n,2},\hat{\cU}_{n,3}),(\hat{\cV}_{n+1},\hat{\cV}_{n,2},\hat{\cV}_{n,3})\right)$, where $\cE_{n,3}$ is defined by  \eqref{eq main error NN generic}. 
\end{Corollary}
}

\section{Numerical results}
\label{se numerics}

{\color{vert}
In this section, we illustrate numerically the theoretical results we have obtained. First, considering the Brownian motion as underlying diffusion, we compare  the schemes described in Section \ref{subse implemented schemes}. Our simulations allow to characterize the discrete time error of each scheme, which corresponds to the theoretical one. 
An empirical study of the computational cost leads us to conclude that using the Crank-Nicholson scheme is optimal when looking for order 2 schemes. Secondly, we test the Euler schemes and the Crank-Nicholson scheme for CIR process as underlying diffusion.
}

\vspace{4mm}
\noindent In the numerical experiments below, the neural network used are  fully connected feedforward network with 2 hidden layers and with $d+10$ number of neurons for each hidden layer. 
A $\mathrm{tanh}$ activation function is used after each hidden layer. The SGD algorithm
\footnote{
We implemented our numerical experiment in Python3 using Tensorflow and adopting  the multi-process techniques to run several tests at the same time. Due to the long running time necessary to characterize the order of convergence,  the code is ran on the LPSM server. We thank the members of the IT team for their support.
}
(ADAM) to perform the optimization is tuned as follows. We set the $batchsize = b$ to train the network and check the convergence of loss function with a test dataset of $batchsize = 2b$  after every 50 training epochs. The learning rate is decreased with a discount factor $\gamma = 0.5$ if the loss  decay is  less than a \textcolor{black}{ given threshold} (set to 0.05 in our case), the training is stopped after the learning rate is less than $10^{-9}$. 
For the Euler schemes, it appears that  choosing $10^{-6}$ yields good results.
However, for high order schemes, we choose a smaller stopping learning rate in order to reduce the impact of the variance of $Y_0$. 
We observed that a small batchsize $b=1000, 10000$ is enough for the Brownian motion case. However, for the CIR process, a small batchsize lead to a bias of $Y_0$ possibly due to the cumulative error from the accuracy of gradient by Monte Carlo simulation.  

Below, we measure the error as the absolute value of the difference between the theoretical solution and average result of $\mathrm{nTest}$ runs as follows:
\[ \epsilon := \left| \frac{1}{\mathrm{nTest}} \sum_{i=1}^{\mathrm{nTest}} \hat{Y}_0^i - Y_0 \right|, \]
and we set $\mathrm{nTest}=10$ in the numerical results below.

\subsection{Comparative study of the schemes}

We study the special case where  $\{\cX_t\}_{0\le t \le T}$ is a drifted Brownian motion, that is
\begin{align}
	\cX_t & = \cX_0 + \mu t + \sigma W_t, \quad 0 \leq t \leq T.
\end{align}
 We use directly the forward diffusion process $\cX$ on the the grid $\Pi$. There is no discretization error for this special case and it corresponds directly to the Euler scheme on $\Pi$.
Since we are in the case where the underlying diffusion is the Brownian motion, we note that $ \cL^{(0)}  \circ \cL^{(\ell)} =  \cL^{(\ell)}  \circ  \cL^{(0)} $ for $ \ell \in \{1, \cdots, d\} $. By Remark 2.1(ii) of \cite{chassagneux2014runge}, one can choose to compute the $Z$-part, the $H$ random weight with this simple form: for $n < N$,
\begin{align*}
	H^n_{q} &= \frac{W_{t_{n+1}} - W_{t_{n,q}} }{c_q h}, \quad 1 < q \le Q+1, \\
	H^n_{q,k} &= \frac{W_{t_{n,k}} - W_{t_{n,q}} }{t_{n,k} - t_{n,q}} = \frac{W_{t_{n,q}} - W_{t_{n,k}} }{(c_k -c_q) h} ,\quad 1<q<k\le  Q+1. 
\end{align*}

\noindent It can be verified that $H^n_{q},  H^n_{q,k} $ satisfy the assumptions $(\mathrm{H}X)_2$, see Proposition 2.3 in \cite{chassagneux2014runge}.

\vspace{2mm}
\noindent We consider the following example borrowed from \cite{hure2020deep}. For $d=10, T=1, t \in [0, T]$, let
\begin{align*}
	\ud\cX_t &= \frac{0.2}{d} \textbf{1}_d \ud t + \frac{1}{\sqrt{d}} \textbf{I}_d \ud W_t, \quad \cX_0 = \textbf{1}_d,  \\
	f(t, x, y, z) &= \left( \cos(\bar{x}) + 0.2 \sin(\bar{x})  \right) e^{\frac{T-t}{2}} - \frac{1}{2}( \sin(\bar{x}) \cos(\bar{x}) e^{T-t})^2 + \frac{1}{2d}\left(y\left(z \cdot \textbf{1}_d \right)\right)^2, \\
	g(x) &= \cos(\bar{x}) , 
\end{align*}
where $ \bar{x} = \sum\limits_{i=1}^d x_i $. The theoretical solution of this BSDE is
$ Y_t = \cos(\bar{X}_t)  e^{\frac{T-t}{2}}$ and $Z_t^i = - \frac{1}{\sqrt{d}}  \sin(\bar{X}_t )  e^{\frac{T-t}{2}}, i = 1, \cdots, d. $\\

\noindent We  compare the numerical results for the five schemes presented in Section \ref{subse implemented schemes} namely the explicit Euler scheme, implicit Euler scheme, Crank-Nicolson scheme, RK-2 scheme and RK-3 scheme. For RK-2 scheme, we set $c_2 = 0.5$ and 
\[ a_{21} = \alpha_{21} = c_2, \quad a_{31} = \alpha_{31} = 1 - \frac{1}{2c_2}, \quad a_{32} = \alpha_{32} = \frac{1}{2c_2},  \]
recall \Cref{subse two stage exp RK}. For RK-3 scheme, we set $c_2 = 0.3, c_3 = 0.7 $ and 
\begin{align*}
	&  a_{21} = \alpha_{21} = c_2,  \quad a_{31} = \alpha_{31} = \frac{c_3(3c_2 - 3c_2^2 - c_3)}{c_2(2-3c_2)}, \quad
	a_{32} = \alpha_{32} = \frac{c_3(c_3 - c_2)}{c_2(2-3c_2)}, \\
	& a_{41} =  \alpha_{41} = \frac{-3c_3 + 6c_2c_3 -3c_2}{6c_2c_3} ,\quad
	a_{42} =  \alpha_{42} = \frac{3c_3 - 2}{6c_2(c_3 - c_2)}, \quad
	a_{43} =  \alpha_{43} = \frac{2-3c_2}{6c_3(c_3 - c_2)},
\end{align*}
recall \Cref{subse three stage exp RK}.

\begin{figure}[H]
\centering
\begin{minipage}[t]{0.42\textwidth}
\centering
\includegraphics[width = 7cm]{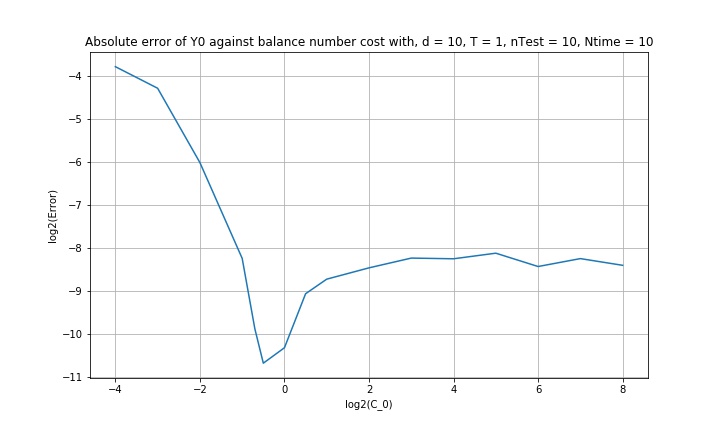} \caption{The error of $Y_0$ for CN scheme based on the balance number}\label{balance}
\end{minipage}
\begin{minipage}[t]{0.1\textwidth}
\centering
\quad
\end{minipage}
\begin{minipage}[t]{0.42\textwidth}
\centering
\includegraphics[width = 7cm]{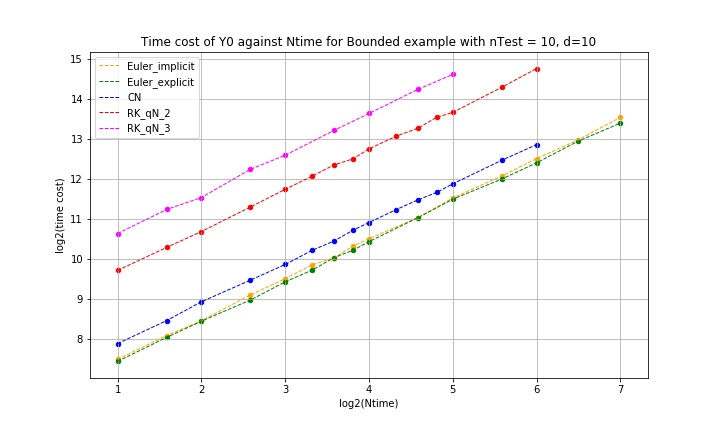} \caption{Time cost against time steps for different schemes} \label{time cost via Ntime}
\end{minipage}
\end{figure}

\vspace{4mm}
\noindent We first quickly discuss the balance number $\mfc$ as it has some influence on the numerical results. On \Cref{balance}, we can see that the error of CN scheme is stable when  $ \log_2(\mfc) > 3 $. However, we can get smaller error when $ -1 < \log_2(\mfc) < 3 $. In the implementation, we choose $\mfc = \frac43$ for CN scheme and $\mfc =  {25 c_q} $ for Runge-Kutta scheme at stage $1<q <Q+1$.

\vspace{2mm}
\noindent \Cref{time cost via Ntime} shows that the computational time is almost linear to the number of time steps for each scheme which is reasonable. The computational time of implicit Euler scheme and explicit scheme are almost the same. The computational time of CN scheme is only slightly larger than Euler scheme since this is still a one-stage scheme though the driver $f$ is computed on both $t_n$ and $t_{n+1}$ at each step $0 \le n \le N-1$. For the Runge-Kutta scheme in the general case: the computational time  is more than $Q$ times the computational time of CN scheme, as expected.\\

\begin{figure}[H]
\centering
\begin{minipage}[t]{0.42\textwidth}
\centering
\includegraphics[width = 7cm]{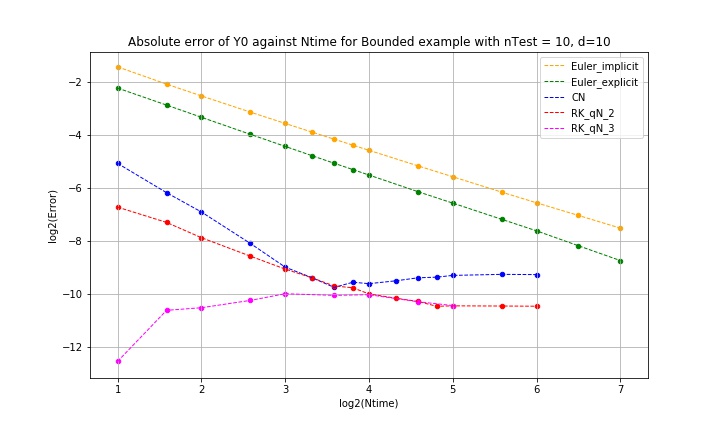} \caption{Error against time steps for different schemes}\label{error via Ntime}
\end{minipage}
\begin{minipage}[t]{0.1\textwidth}
\centering
\quad
\end{minipage}
\begin{minipage}[t]{0.42\textwidth}
\centering
\includegraphics[width = 7cm]{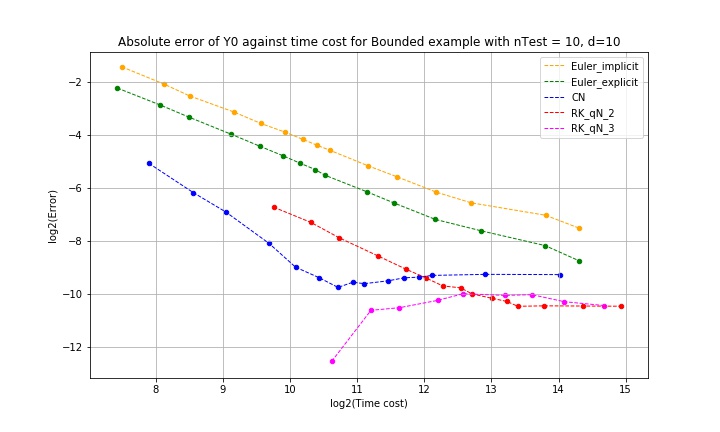} \caption{Error against time cost for different schemes} \label{error via time cost}
\end{minipage}
\end{figure}

In \Cref{error via Ntime}, we compare the convergence rate of the five schemes mentioned above. We verify that the implicit Euler scheme and explicit Euler scheme are almost order 1. The CN scheme is almost order 2. The convergence rate of RK-2 scheme is slightly less than CN scheme, but the error is smaller. The RK-3 scheme converges too fast in terms of discrete time error for us to be able to observe any convergence order.  \\

In \Cref{error via time cost},  we plotted the error w.r.t. the time cost for the five schemes mentioned above. We see that the Euler schemes are slower to reach a small error. The RK-3 scheme has a small error but it is computationally demanding even if the number of time steps $N$ is small. As expected, the CN scheme is faster than RK-2 scheme. In conclusion, if we want an error smaller than $ 0.01 \approx 2^{-6.64} $, CN scheme seems to be the best scheme to use.

\subsection{An application with Cox–Ingersoll–Ross process}
In this section, we test the Cox–Ingersoll–Ross \cite{cozma2018} (CIR for short) process as follows,
\[ d \cX_t = a(b - \cX_t) dt + \sigma \sqrt{\cX_t} dW_t, \]
where $\sqrt{\cX_t} =  {\rm diag}( \sqrt{\cX_t^1}, \cdots,  \sqrt{\cX_t^d}) $, and with the following conditions to ensure that $ (\cX_t)_{t>0} $ is always positive:
\[ a>0,\; b >0, \; 2ab \ge \sigma^2,\; \cX_0 \ge 0. \]

\noindent We implemented the Ninomiya-Victoir scheme \cite{ninomiya2008weak}: The weak error is of order $2$, recall \eqref{eq approx exp y}. 


\noindent We test the BSDE with the same solution of the previous subsection $  Y_t = u(t, X_t) = \cos(\bar{X}_t)  e^{\frac{T-t}{2}} $, $Z_t^i = - \sqrt{\frac{X_t^i}{d}} \sin(\bar{X}_t )  e^{\frac{T-t}{2}},  i = 1, \cdots, d$, recalling $\bar{x} := \sum\limits_{i=1}^{d}x_i$, and keep the terminal function $g$ the same with Brownian Motion case, but with the forward diffusion process is CIR process:
\begin{align*}
	\ud\cX_t &= \frac{1}{5d} (3 -  \cX_t) \ud t + \frac{1}{\sqrt{d}} \sqrt{\cX_t}   \ud W_t, \quad \cX_0 = 10 \textbf{1}_d.
\end{align*}

Then,  we set the driver $f$:
\begin{align*}
	\tilde{f} (t, x)  &= \left(  \frac{1}{2}  \cos(\bar{x}) (1+ \frac{\bar{x}}{d} ) +  \sin(\bar{x}) (\frac35 - \frac{\bar{x}}{5d})  \right) e^{\frac{T-t}{2}}, \\
	f(t, x, y) &= \tilde{f}(t, x)  + \frac{1}{5} \left(\sin(\bar{X}_t ) e^{\frac{T-t}{2}} \right)^2 ( y^2 - \cos(\bar{x})^2 e^{T-t}  ) , 
\end{align*}
where $ \bar{x} = \sum\limits_{i=1}^d x_i $.

 We only compare Crank-Nicolson scheme with implicit Euler scheme in this section. Setting $d=10, T=1, \mfc = 1$. We test the implicit Euler scheme and Crank-Nicolson scheme   with $b= 50000$. \\
 
\begin{figure}[H]
\centering
\begin{minipage}[t]{0.42\textwidth}
\centering
\includegraphics[width = 7cm]{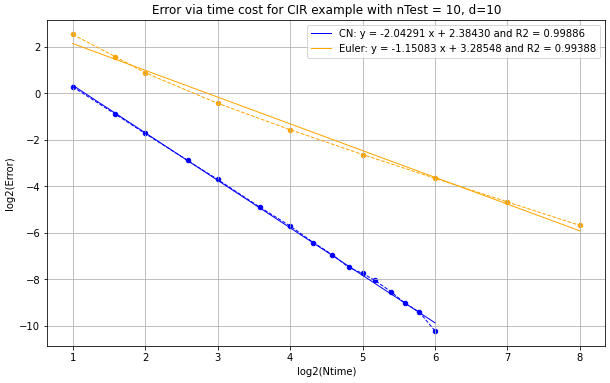} \caption{Error against time steps for Euler and CN schemes in CIR case}\label{CIR y error via Ntime}
\end{minipage}
\begin{minipage}[t]{0.1\textwidth}
\centering
\quad
\end{minipage}
\begin{minipage}[t]{0.42\textwidth}
\centering
\includegraphics[width = 7cm]{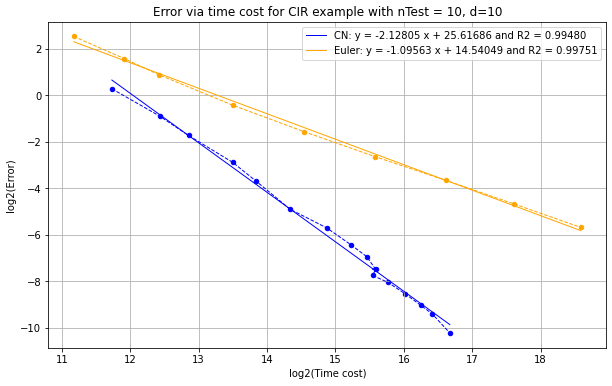} \caption{$Y_0$ against time steps for for Euler and CN schemes in CIR case} \label{CIR y error via time}
\end{minipage}
\end{figure}

 In \Cref{CIR y error via Ntime} and \Cref{CIR y error via time}, the implicit Euler scheme is almost order 1, however it can only achieve an error around $2^{-6}$ even $N=256$, which spent more than $2^{18}$ seconds (nearly 5 days).  And the Crank-Nicolson scheme is almost order 2, the error can achieve $ 2^{-10} $ when $N=64$, and the error can be less than $ 2^{-6} $ with $ 2^{15} $ seconds (about 9 hours). \\

\section{Appendix}

\subsection{Proof of Proposition \ref{pr generic stab for cond exp estim}} \label{subse proof of prop}
Recalling (\ref{eq de interm stage Y})-(\ref{eq de interm stage Z}) and $ \upsilon^n_q = \EFp{t_{n,q}}{ |H_q^n|^2}  $,  observe that for $1 < q \le Q+1$, it holds
\begin{align}\label{Yij RK qN}
	Y_{n,q} =  Y_{n+1} +  h \sum_{k=1}^{q} a_{qk} f(X_{n,k}, Y_{n,k}, Z_{n,k})  - (Z_{n, q} + A_{n,q} )\Delta W_{n, q}  + \Delta M_{n,q},
\end{align}
where $  \EFp{t_{n,q}}{ \Delta M_{n, q} } =  \EFp{t_{n, q}}{\Delta M_{n, q} \Delta W_{n, q}  } =0,  \EFp{t_{n, q}}{ |\Delta M_{n, q} |^2} < \infty $, 
and
\begin{align}\label{definition A_ij RK qN}
	A_{n, q} &= \EFp{t_{n, q}}{  H_q^n \left( Y_{n+1} +  h \sum_{k=1}^{q-1} a_{qk} f(X_{n,k}, Y_{n,k}, Z_{n,k}) \right)} - Z_{n, q} \nonumber \\
	&= \EFp{t_{n, q}}{  \sum_{k=1}^{q-1}  \left( a_{qk}  H_q^n -\alpha_{qk}  H_{q,k}^n \right) h f(X_{n,k}, Y_{n,k}, Z_{n,k})   }  .
\end{align}
And for the perturbed scheme defined in (\ref{eq de interm stage Y perturbed})-(\ref{eq de interm stage Z perturbed}), observe that
\begin{align}\label{Yij RK qN perturbed}
	\tilde{Y}_{n,q} =  \tilde{Y}_{n+1} +  h \sum_{k=1}^{q} a_{qk} f(X_{n,k}, \tilde{Y}_{n,k}, \tilde{Z}_{n,k}) + \zeta^y_{n,q} - (\tilde{Z}_{n, q} + \tilde{A}_{n,q} - \zeta^z_{n,q})\Delta W_{n, q} + \Delta \tilde{M}_{n,q},
\end{align}
where $  \EFp{t_{n,q}}{ \Delta \tilde{M}_{n, q} } =  \EFp{t_{n, q}}{ \Delta \tilde{M}_{n, q} \Delta W_{n, q} } =0,  \EFp{t_{n, q}}{ |\Delta \tilde{M}_{n, q} |^2} < \infty $, and 
\begin{align}\label{definition A_ij RK qN}
	\tilde{A}_{n, q} = \EFp{t_{n, q}}{  \sum_{k=1}^{q-1}  \left( a_{qk}  H_q^n -\alpha_{qk}  H_{q,k}^n \right) h f(X_{n,k}, \tilde{Y}_{n,k}, \tilde{Z}_{n,k})   }  .
\end{align}
 Set $\delta Y_{n, q} = Y_{n,q} - \tilde{Y}_{n,q},  \delta Z_{n, q} = Z_{n, q} - \tilde{Z}_{n,q}, \delta A_{n, q} = A_{n, q} - \tilde{A}_{n, q}, \delta f_{n, q} = f(X_{n, q}, Y_{n, q}, Z_{n, q}) -   f(X_{n, q}, \tilde{Y}_{n, q}, \tilde{Z}_{n, q}) $, $ \delta \Delta M_{n, q} =  \Delta M_{n, q} - \Delta \tilde{M}_{n, q}$ for all $ 0\le n \le N-1, 1 <  q\le Q+1 $. From \cref{Yij RK qN} and \cref{Yij RK qN perturbed}, we get
\begin{align}\label{difference of two equations}
 	\delta Y_{n, q} +  (\delta Z_{n, q} + \zeta_{n, q}^{z} ) \Delta W_{n, q} + \delta \Delta M_{n, q} = \delta Y_{n+1} + h \sum_{k=1}^q a_{qk} \delta f_{n, k} - \zeta_{n, q}^{y}   - \delta A_{n, q}  \Delta W_{n, q} .
 \end{align}
\textit{Step 1: For $0 \le n \le N-1$, control of $ \esp{|\delta Y_n|^2}$ by  the term $ h \sum\limits_{k=n}^{N-1} \sum\limits_{q=2}^{Q+1} \esp{ |\delta Z_{k,q}|^2 }$.}\\
Squaring both sides of \eqref{difference of two equations}, taking conditional expectation and  using Young’s inequality, we obtain 
 \begin{align}\label{square both side}
 	|\delta Y_{n, q} |^2 + c_q h |\delta Z_{n, q} + \zeta_{n, q}^{z} |^2 & \le  \left(1 + \frac{h}{C} \right) \EFp{t_{n, q}}{ \left|\delta Y_{n+1} + h \sum_{k=1}^j a_{qk} \delta f_{n, k}  \right|^2  } \nonumber \\ 
 	&+ 2\left(1 + \frac{C}{h} \right) \left( \left| \zeta_{n, q}^{y} \right|^2 +  c_qh \left| \delta A_{n, q} \right|^2 \right).
 \end{align}
 Recalling (\ref{eq conditions on H 1})-(\ref{eq conditions on H 2}), and denote $\bar{a} = \max\limits_{q,k}\{ |a_{qk}|, |\alpha_{qk}| \} $, thus have
 \begin{align} \label{control A_ij}
 	\left| \delta A_{n, q}  \right|^2 &=  h^2 \left|\EFp{t_{n, q}}{ \sum\limits_{k=1}^{q-1} (a_{qk} H_q^n- \alpha_{qk}  H_{q,k}^n ) \delta f_{n,k} } \right|^2  \nonumber  \\
 	&  \le q  h^2 \sum\limits_{k=1}^{q-1}  \EFp{t_{n, q}}{ | a_{qk} H_q^n- \alpha_{qk}  H_{q,k}^n  |^2 }  \EFp{t_{n, q}}{ |\delta f_{n,k} |^2 } \nonumber \\
 	& \le  q  h^2 \sum\limits_{k=1}^{q-1}  \EFp{t_{n, q}}{ 2 a_{qk}^2|  H_q^n |^2 +2 \alpha_{qk} ^2 |  H_{q,k}^n  |^2 }  \EFp{t_{n, q}}{ 2 K^2( |\delta Y_{n, k} |^2 + |\delta Z_{n,k} |^2) } \nonumber \\
 	& \le 8q^2\bar{a}^2  K^2 \Lambda  h \sum\limits_{k=1}^{q-1}  \EFp{t_{n, q}}{ |\delta Y_{n,k} |^2 + |\delta Z_{n,k} |^2 }  \nonumber \\
 	& \le C  h \sum\limits_{k=1}^{q-1}  \EFp{t_{n, q}}{ |\delta Y_{n,k} |^2 + |\delta Z_{n,k} |^2 } ,
 \end{align}
 and
 \begin{align}\label{control delta Z_ij and zeta_ij}
 	|\delta Z_{n, q} + \zeta_{n, q}^{z} |^2 \ge (1 - \eta )  |\delta Z_{n, q}|^2 + (1-\frac{1}{\eta}) |\zeta_{n, q}^{z}|^2 \ge (1 - \eta )  |\delta Z_{n, q}|^2 -\frac{1}{\eta} |\zeta_{n, q}^{ z}|^2,
 \end{align}
  for $ \eta >0 $. Using the Lipschitz regularity of $f$, Young's inequality and Jensen's inequality, for any $\epsilon > 0$, we obtain
 \begin{align*}
  & \quad\,	\left|\delta Y_{n+1} + h \sum_{k=1}^q a_{qk} \delta f_{n, k}  \right|^2 \\
   	& \le  \left( |\delta Y_{n+1}| +  h K \sum_{k=1}^q |a_{qk}| ( |\delta Y_{n, k} | + |\delta Z_{n, k} | ) \right)^2 \nonumber \\
 	& \le  \left( \left(1 + C h \right) |\delta Y_{n+1}|  + C h   |\delta Z_{n+1} | + Ch \sum_{k=2}^q  ( |\delta Y_{n, k} | + |\delta Z_{n, k} | ) \right)^2 \nonumber \\
 	& \le \left(1 + \frac{h}{\epsilon} \right)   \left( \left(1 + C h \right) |\delta Y_{n+1}|^2 +   Ch  \sum_{k=2}^q  |\delta Y_{n, k} |^2 \right) + C \left(h + \epsilon \right) h  \sum_{k=1}^q  |\delta Z_{n, k} |^2  .
 \end{align*}
Choosing $h$ small enough and $\epsilon$ such that $ C (h+\epsilon) \le \frac{c_q}{2 } $, we obtain
 \begin{align}\label{control delta Y_ij}
 	\left|\delta Y_{n+1} +  h \sum_{k=1}^q a_{qk} \delta f_{n, k}  \right|^2  & \le \left(1 + Ch \right) |\delta Y_{n+1}|^2 + Ch   \sum_{k=2}^q  |\delta Y_{n, k} |^2 + \frac{c_q h}{2 }   \sum_{k=1}^q  |\delta Z_{n, k} |^2  .
 \end{align}
 Thus choosing $\eta = \frac{1}{4} $ in (\ref{control delta Z_ij and zeta_ij}), and using (\ref{control delta Y_ij}), (\ref{control A_ij})  into (\ref{square both side}), then for $h$ small enough,  we get
  \begin{align} \label{upper bound delta Y_ij}
  & \quad\,	|\delta Y_{n, q} |^2  \le |\delta Y_{n, q} |^2 + \frac{c_q h}{4 }  |\delta Z_{n, q} |^2  \nonumber \\
 	 &\le \left(1 + Ch \right) \EFp{t_{n, q}}{ |\delta Y_{n+1}|^2} + Ch   \sum_{k=1}^{q - 1}  \EFp{t_{n, q}}{ |\delta Y_{n, k} |^2 + |\delta Z_{n, k} |^2 } + \EFp{t_{n, q}}{  \frac{C}{h} | \zeta_{n, q}^{y} |^2  +  Ch \left| \zeta_{n, q}^{z} \right|^2 }  
 	 \nonumber \\
 	 &\le \left(1 + Ch \right) \EFp{t_{n, q}}{ |\delta Y_{n+1}|^2}  + Ch   \sum_{k=1}^{q - 1}  \EFp{t_{n, q}}{  |\delta Z_{n, k} |^2 } +C \sum_{k=2}^{q} \EFp{t_{n, q}}{ \frac{1}{h} | \zeta_{n,k}^{y} |^2  +  h \left| \zeta_{n, k}^{z} \right|^2 } .
 \end{align}
  Using the discrete version of Grönwall's lemma, we even eventually conclude:
	\begin{align} \label{upper bound delta Y prox RK}
		\esp{|\delta Y_n|^2}  &\le C\left( \esp{ |\delta Y_N|^2 + |\delta Z_N|^2  }   +  h\sum_{k=n}^{N-1}  \sum_{q=2}^{Q+1} \esp{    |\delta Z_{k, q} |^2  +   \frac{1}{h^2} | \zeta_{k, q}^{y} |^2  +  \left| \zeta_{k, q}^{z} \right|^2   } \right) .
	\end{align}
\noindent
 \textit{Step 2: Control of $ h \sum\limits_{k=n}^{N-1}  \sum\limits_{q=2}^{Q+1}\esp{ |\delta Z_{k, q}|^2} $.} \\ 
 Using the Cauchy-Schwarz inequality and the Lipschitz regularity of $f$, we get
 \begin{align}\label{upper bound 1 delta Zij}
 	 h \left|\EFp{t_{n, q}}{ h H_{q,k}^n \delta f_{n,k} } \right|^2  	& \le h^3 \EFp{t_{n, q}}{ |H_{q,k}^n|^2}  \EFp{t_{n, q}}{ |\delta f_{n,k} |^2}    \le 2dK^2\Lambda h^2 \EFp{t_{n, q}}{ |\delta Y_{n,k} |^2 + |\delta Z_{n,k} |^2 }
 \end{align}
 and 
 \begin{align}\label{upper bound 2 delta Zij}
 \EFp{t_n}{ h \left| \EFp{t_{n, q}}{  H_q^n (\delta Y_{n+1}  -   \EFp{t_{n}}{\delta Y_{n+1}} )  } \right|^2 }
 & \le   \EFp{t_n}{ h \EFp{t_{n, q}}{ |H_q^n|^2 } \EFp{t_{n, q}}{ |\delta Y_{n+1}  -   \EFp{t_{n}}{\delta Y_{n+1}} |^2}   }\nonumber \\
 & \le \Lambda \EFp{t_n}{  |\delta Y_{n+1}  -   \EFp{t_{n}}{\delta Y_{n+1}} |^2}   \nonumber \\
 &\le \Lambda \left( \EFp{t_{n}}{ |\delta Y_{n + 1} |^2} - \EFp{t_{n}}{\delta Y_{n+1}} ^2 \right) .
 \end{align}
 Taking $q = Q+1$ in \eqref{difference of two equations}, using Young's inequality, and using the subscript {$ \{n\} \equiv \{n, Q+1\}$}, we get
 \begin{align}\label{upper bound 3 delta Zij}
 	 & \quad\, \EFp{t_{n}}{ |\delta Y_{n + 1} |^2} - \EFp{t_{n}}{\delta Y_{n+1}}^2   \nonumber \\
 	& =  \EFp{t_{n}}{ |\delta Y_{n + 1} |^2} -  \EFp{t_n}{ \delta Y_{n} + \zeta_{n}^{y} - h\sum_{k=1}^{Q+1} a_{qk}  \EFp{t_{n, q}}{ \delta f_{n, k} }  }^2 \nonumber \\
 	& \le   \EFp{t_{n}}{ |\delta Y_{n + 1} |^2} -  |\delta Y_{n} |^2  + C\left(  |\delta Y_{n} | |\zeta_{n}^{y}|  + h ( |\delta Y_{n} | +  |\zeta_{n}^{y} |  )\sum_{k=1}^{Q+1}  \EFp{t_{i}}{ |\delta f_{n, k}| }  \right)   \nonumber \\
 	& \le   \EFp{t_{n}}{ |\delta Y_{n + 1} |^2} -  |\delta Y_{n} |^2  + C\left( |\delta Y_{n} | |\zeta_{n}^{y}|  + h ( |\delta Y_{n} | +  |\zeta_{n}^{y} |  ) \sum_{k=1}^{Q+1}  \EFp{t_{n, q}}{  |\delta Y_{n,k} | + |\delta Z_{n,k} | } \right)   \nonumber \\
 	& \le   \EFp{t_{n}}{ |\delta Y_{n + 1} |^2} -  |\delta Y_{n} |^2   + Ch  \sum_{k=1}^{Q+1}  \EFp{t_{n}}{ \frac{1}{\epsilon} |\delta Y_{n,k} |^2 + \epsilon |\delta Z_{n,k} |^2 }   + \frac{C}{h}  \EFp{t_{n}}{ |\zeta_{n}^{y}|^2 }, 
 \end{align}
 for any $\epsilon>0$.
 Using Jensen's inequality, and the inequalities (\ref{upper bound 1 delta Zij}), (\ref{upper bound 2 delta Zij}), (\ref{upper bound 3 delta Zij}), we obtain
 \begin{align}
 	& \quad\, h \EFp{t_n}{ |\delta Z_{n, q}|^2 } = h \EFp{t_n}{ \left| \mathbb{E}_{t_{n, q}}\Big[ H_{q}^n \delta Y_{n+1}  + h \sum\limits_{k=1}^{q-1} \alpha_{qk} H_{q,k }^n \delta f_{n,k}  - \zeta_{n, q}^{z}  \Big] \right|^2 } \nonumber \\
 	 &\le C \EFp{t_n}{  h \left| \EFp{t_{n, q}}{ H_q^n (\delta Y_{n+1}  -  \EFp{t_{n}}{\delta Y_{n+1}} )  } \right|^2+ \sum_{k=1}^{q-1}  \alpha_{qk}^2 h \left|\EFp{t_{n, q}}{ h H_{q,k}^n \delta f_{n,k} } \right|^2 + h|\zeta_{n,q}^{z}|^2  } \nonumber \\
 	 & \le C \left( \EFp{t_{n}}{ |\delta Y_{n + 1} |^2} - |\EFp{t_{n}}{\delta Y_{n+1}} |^2   + h^2  \sum_{k=1}^{q-1}  \EFp{t_{n}}{ |\delta Y_{n,k} |^2 +  |\delta Z_{n,k} |^2 }    + h  \EFp{t_{n}}{ |\zeta_{n, q}^{z}|^2 }  \right) \nonumber \\
 	 & \le C \left( \EFp{t_{n}}{ |\delta Y_{n + 1} |^2} -  |\delta Y_{n} |^2   + Ch  \sum_{k=1}^{Q+1}  \EFp{t_{n}}{ \frac{1}{\epsilon} |\delta Y_{n,k} |^2 + \epsilon |\delta Z_{n,k} |^2 }     +   \EFp{t_{n}}{\frac{C}{h} |\zeta_{n}^{y}|^2 +    h |\zeta_{n,q}^{z}|^2 }  \right) 
 \end{align}
Summing over $q$ and $n$, setting $ \epsilon = \frac{1}{2(Q+2)C} $, and note that the subscript  {$ \{n, 1\} \equiv \{n+1\} \equiv \{n+1, Q+1\} $}, then for $h$ small enough,
\begin{align*}
	& \quad\, h \sum_{k=n}^{N-1} \sum_{q=2}^{Q+1} \EFp{}{ |\delta Z_{k, q}|^2} \nonumber \\
	  &\le C \sum_{k=n}^{N-1}  \sum_{q=2}^{Q+1}  \EFp{}{  |\delta Y_{k + 1} |^2 -  |\delta Y_{k} |^2    + h  \sum_{l=1}^{Q+1} ( \frac{1}{\epsilon}|\delta Y_{k,l} |^2 + \epsilon |\delta Z_{k,l} |^2  ) + \frac{1}{h} |\zeta_{k}^{y} |^2  + h|\zeta_{k,q}^{z}|^2  } \nonumber \\
	& = C \sum_{k=n}^{N-1}   \EFp{}{  Q|\delta Y_{k + 1} |^2 -  Q|\delta Y_{k} |^2    + \frac{Qh}{\epsilon}  \sum_{l=1}^{Q+1} |\delta Y_{k,l} |^2   + \frac{Q}{h} |\zeta_{k}^{y} |^2  + h  \sum_{q=2}^{Q+1}  |\zeta_{k, q}^{z}|^2  } \nonumber \\
	&\quad + \frac{Qh}{Q+2}  \sum_{k=n}^{N-1}  \sum_{l=2}^{Q+1} \esp{   |\delta Z_{k,l} |^2 }  
\end{align*}
 so that
\begin{align}\label{upper bound Zij}
	& \quad\, h \sum_{k=n}^{N-1} \sum_{q=2}^{Q+1} \EFp{}{ |\delta Z_{k,q}|^2} \nonumber \\
	& \le  C \sum_{k=n}^{N-1}   \EFp{}{  |\delta Y_{k + 1} |^2 -  |\delta Y_{k} |^2     + h \sum_{q=1}^{Q+1} |\delta Y_{k,q} |^2   + \frac{1}{h} |\zeta_{k}^{y} |^2  + h  \sum_{q=2}^{Q+1}  |\zeta_{k, q}^{z}|^2  }  \nonumber \\
	& \le C\esp{ |\delta Y_N|^2 + h |\delta Z_N |^2 +h  \sum_{k=n}^{N-1}   |\delta Y_{k} |^2  + \sum_{k=n}^{N-1} \sum_{q=2}^{Q+1}  \left( \frac{1}{h} |\zeta_{k, q}^{y} |^2  + h|\zeta_{k, q}^{z}|^2 \right)  }
\end{align}
where we used (\ref{upper bound delta Y_ij}) and the subscript  {$ \{n, 1\} \equiv \{n+1\} \equiv \{n+1, Q+1\}, 0\le n\le N-1 $} again in the last line.\\
\noindent
 \textit{Step 3: Control of the term $ \esp{|\delta Y_n|^2} , \; 0 \le n \le N-1 $.} \\ 
 Combining the inequality (\ref{upper bound Zij}) with (\ref{upper bound delta Y prox RK}), we get
\begin{align}
	\esp{|\delta Y_n|^2}  &\le C \esp{ |\delta Y_N|^2 + h |\delta Z_N|^2 + h  \sum\limits_{k=n}^{N-1}  |\delta Y_k |^2  +  h\sum_{k=n}^{N-1} \sum_{q=2}^{Q+1}  \left( \frac{1}{h^2} |\zeta_{k, q}^{y} |^2  + |\zeta_{k, q}^{z}|^2 \right)   } .
\end{align}
\textit{Step 4: Control on $\max\limits_{n\le k \le N-1} \esp{|\delta Y_k|^2} + h \sum\limits_{k=n}^{N-1} \esp{ |\delta Z_k|^2 } $. }\\
Set $ \delta_n : = \sum\limits_{k=n}^{N-1} \esp{  |\delta Y_k |^2 } $, and 
\begin{align}
	\theta_n := \esp{ |\delta Y_N|^2 + h |\delta Z_N|^2} +   h\sum_{k=n}^{N-1}\sum_{q=2}^{Q+1} \esp{\frac{1}{h^2} |\zeta_{k, q}^{y} |^2  + |\zeta_{k, q}^{z}|^2  } .
\end{align}
It holds
\begin{align} \label{control delta_n}
	\delta_n - \delta_{n+1} \le Ch \delta_n + C \theta_n .
\end{align}
Using the discrete version of Grönwall's lemma, and noting that $ \delta_{N-1}  \le \theta_n$, $\theta_k \le \theta_n$ for $ k \ge n$, we obtain
\begin{align} \label{upper bound delta_n}
	\delta_n \le C \left( \delta_{N-1} + \sum\limits_{k=n}^{N-1} \theta_k e^{C(N-k-1)h} \right) \le C \theta_n  \frac{1}{e^{Ch}-1}.
\end{align}
This last inequality combined with \eqref{control delta_n} leads to 
\begin{align*}
	\esp{  |\delta Y_n |^2 } = \delta_n - \delta_{n+1} \le C \theta_n.
\end{align*}
For the $Z-$part, the proof is concluded inserting  \eqref{upper bound delta_n} into \eqref{upper bound Zij} with $n=0$ in this equation. 

\bibliographystyle{abbrv}
\bibliography{shortbib}

\end{document}